\documentclass[english, 11pt]{amsart} 
\usepackage{std_macros_new}
\newcommand{\SU}{\textup{SU}}

\title{The trace formula and the existence of PEL-type Abelian varieties modulo $p$}

\author{Arno Kret}

\let\mathscr\mathcal
\begin{document}

\begin{abstract}
We show, using the trace formula, that any Newton stratum of a Shimura variety of PEL-type of types (A) and (C) is non-empty at the primes of good reduction. Furthermore we prove \emph{conditionally} the non-emptiness for Shimura data associated to odd Spin groups.

Our results are conditional on Rapoport-Langlands conjecture and Arthur's conjectures on the discrete spectrum. Both these results have been announced by Arthur and Kisin in significant cases.  
\end{abstract}

\maketitle

\tableofcontents


\section*{Introduction} 

Consider a Shimura variety of PEL-type and reduce it modulo a prime $p$ where this variety has good reduction. Then the variety parametrizes Abelian varieties in characteristic $p$ with certain additional PEL-type structures. To each such Abelian variety one may associate its Dieudonn\'e isocrystal. The PEL-structures on the Abelian variety give additional $G$-structure on the isocrystal, and as such the isocrystals lie in the category of ``isocrystals with additional structures'' (Kottwitz \cite{MR809866}). We look at these objects modulo isomorphism. When given a Shimura variety of PEL-type modulo $p$, not every isocrystal with additional $G$-structure arises from a geometric point on this variety. In fact, there are only a finite number of possible isocrystals; since the work of Rapoport-Richarz and Kottwitz \cite{MR1411570, MR1485921} it is known that they all lie in a certain, explicit finite set $B(G_{\qp}, \mu)$ of `admissible' isocrystals. However, they did not show that $B(G_{\qp}, \mu)$ is \emph{exactly} the set of possible isocrystals, \ie that for every element $b \in B(G_{\qp}, \mu)$ there exists some Abelian variety in characteristic $p$ with additional PEL-type structures whose rational Dieudonn\'e module is equal to $b$. Recently Viehmann and Wedhorn  \cite{viehmann} have proved through geometric means that this is indeed the case if the group of the Shimura datum is of type (A) or (C), and also recently Sug Woo Shin announced a proof. 

The article is organized as follows. In Chapter 1 we introduce some generalities on truncated traces on smooth representations of certain $p$-adic groups. These results are the main local ingredients that go into the global proofs of the later chapters.

In Chapter 2 we give a new proof of the non-emptiness result using automorphic forms and the trace formula for Shimura varieties of type (A).


In Chapter 3 we show non-emptiness results for larger classes of Shimura varieties that are not necessarily of PEL type. However, our arguments are \emph{conditional on two deep conjectures}: (C1) A suitable form of the Kottwitz formula for the Shimura datum $(G, X)$ at the prime of reduction $p$, (C2) A suitable form of the conjecture of Arthur on the discrete spectrum of $G$. Both (C1) and (C2) have been announced in significant cases, by Vasiu, Kisin (for (C1)) and Arthur (for (C2)). We thus have good hope that our results will become unconditional soon in the respective significant cases. In fact, we do not need a condition as strong as (C2); we only need a much weaker variant of it. The condition is too technical to state in this introduction, see Hypothesis~\ref{Ghyp} in the text. 

In Chapters 4 and 5 we use Arthur's recent book~\cite{arthurbook} to prove Hypothesis~\ref{Ghyp} for the symplectic groups of similitudes and the odd spinor groups of similitudes. We should mention that Arthur's book is currently conditional as well. Arthur's proof relies on the stabilization of the twisted trace formula for the general linear group. Thus, once this stabilization is established, our nonemptiness result becomes dependent only on (C1).

For PEL varieties of type C, conjecture (C1) is proved by Kottwitz. Thus in this case our non-emptiness argument depends only on the stabilization of the twisted trace formula.

Let us explain our method of proof for the unitary group of similitudes, as our method is similar for the other groups. The formula of Kottwitz for the number of points on a Shimura variety modulo $p$ can be restricted to count the number of points in any given Newton stratum. Thus, the problem is to show that this formula of Kottwitz does not vanish after restriction. The formula of Kottwitz is in terms of (twisted) orbital integrals on the group $G$, and he rewrites the formula in terms of stable orbital integrals on certain endoscopic groups of $G$. The advantage of this stable expression is that it may be compared to the geometric side of the stable trace formula. By doing so, one will get a certain sum over the endoscopic groups of certain truncated, transferred Hecke operators acting on automorphic representations of these endoscopic groups. (The truncation is defined by the element of $B(G_{\qp}, \mu)$.) A general objective is to try and work out this expression; one will then get a rather precise description (of the alternating sum) of the cohomology of the Newton strata. In recent work we have done this for certain simple Shimura varieties, cf. \cite{kret1}, \cite{kret2}. In this article we have aimed at a simpler goal: We do not try to describe the cohomology of the Newton stratum defined by $b \in B(G_{\qp}, \mu)$, we only want to show that this cohomology does not vanish, so that the corresponding Newton stratum must be non-empty. This means that we pick one very particular Hecke operator $f^p$ and carry out the computation sketched above only for this particular Hecke operator. We choose our Hecke operator with care, so that all the proper endoscopy vanishes and that in the end, after applying a simple version of the trace formula, we arrive at a sum of certain $b$-truncated traces on cuspidal automorphic representations of the quasi-split inner form $G^*$ of the group $G$ (Equation~\eqref{eqn:Automorphic_Sum_A}):
\begin{equation}\label{eqn:Intro_Trace_on_Automorphic_reps}
\sum_\Pi m(\Pi) \cdot \Tr( ( f^p)^{G^*}, \Pi^p) \Tr(\chi_b^{G(\qp)} f_\alpha, \Pi_p). 
\end{equation}
In fact, the function $f^p$ will be chosen so that there is \emph{at least one} contributing representation $\Pi_0$, and so that for any other hypothetical $\Pi$ contributing to Equation~\eqref{eqn:Intro_Trace_on_Automorphic_reps}, the quotient
\begin{equation}\label{eqn:Intro_Same_Signs}
\frac {m(\Pi) \cdot \Tr( (f^p)^{G^*}, \Pi^p) \Tr(\chi_b^{G(\qp)} f_\alpha, \Pi_p) }{m(\Pi_0) \cdot \Tr( (f^p)^{G^*}, \Pi_{0}^p) \Tr(\chi_b^{G(\qp)} f_\alpha, \Pi_{0,p})},
\end{equation}
is a positive real number in case $\alpha$ is sufficiently divisible. Then, the sum in Equation~\eqref{eqn:Intro_Trace_on_Automorphic_reps} is non-zero (for $\alpha$ sufficiently divisible). Thus the formula of Kottwitz does not vanish as well, and this means that the corresponding Newton stratum is nonempty. An important step in the argument is to show that $\Pi_0$ exists. In particular one has to find a local representation $\Pi_{0, p}$ at $p$ such that for this local representation we have $\Tr(\chi_b^{G(\qp)} f_\alpha, \Pi_{0, p}) \neq 0$. We find a set of such representations with positive Plancherel measure. General theory of automorphic forms then allow us to find a global automorphic representation $\Pi_0$ with $\Pi_{0,p}$ lying in our Plancherel set.

\bigskip

A convention: We often mention integers $\alpha$ which are `sufficiently divisible' such that a certain statement ($\star$) holds. With this, we mean to say, by definition, that there exists an integer $M$ such that whenever $M$ divides $\alpha$, the statement ($\star$) is true for the integer $\alpha$.

\smallskip
When ambiguity is not possible we confuse an algebraic group $G$ defined over a field $F$ with its set of $F$-rational points.
\smallskip

Let $x \in \R$ be a real number, then $\lfloor x \rfloor$ (\emph{floor function}) (resp. $\lceil x \rceil$, \emph{ceiling function}) denotes the unique integer in the real interval $(x-1, x]$ (resp. $[x, x+1)$). 

For us, a \emph{composition} of a positive number $n$ is an \emph{ordered} list of positive integers $(n_a)$ such that $\sum_{a=1}^k n_a = n$. A \emph{partition} is a non-ordered list of non-negative integers summing up to $n$. This convention follows (for example) the books of Stanley.

\bigskip

\noindent {\textbf{Acknowledgements: } I would like to thank my thesis advisers Laurent Clozel and Laurent Fargues. Especially Laurent Clozel who helped me putting everything together and encouraged to me that, after the local combinatorial problems have been resolved, the global part of the argument should be ``de la cuisine''. I also wish to thank Bas Edixhoven, Guy Henniart, Erez Lapid, Giancarlo Lucchini, Judith Ludwig, Ben Moonen, Frans Oort, Peter Scholze, Sug Woo Shin, Eva Viehmann, Torsten Wedhorn, and Chao Zhang.

\section{Truncated traces}
We use truncated traces to study the Newton strata of Shimura varieties. 

\subsection{Isocrystals}\label{sect:isocrystals}
We start this preliminary section with some notations. Let $p$ be a prime number and let $F$ be a finite extension of $\qp$. Let $\cO_F$ be the ring of integers of $F$, let $\varpi_F \in \cO_F$ be a prime element. We write $\Fq$ for the residue field of $\cO_F$, and the number $q$ is by definition its cardinality. We fix an algebraic closure $\lqp$ of $F$, and we let $F_\alpha$ be the unramified extension of $F$ of degree $\alpha$ in $\lqp$. Let $G$ be a smooth reductive group over $\cO_F$ (then $G_F$ is an unramified group \cite{MR546588}). We fix a minimal parabolic subgroup $P_0$ of $G$, and we standardize the parabolic subgroups of $G$ with respect to $P_0$. We write $T \subset P_0$ for the Levi component of $P_0$ and $N_0$ for the unipotent part, so that we have $P_0 = TN_0$. We call a parabolic subgroup $P$ of $G$ \emph{standard} if it contains $P_0$, and we write $P = MN$ for its standard Levi decomposition. We write $K$ for the hyperspecial subgroup $G(\cO_F) \subset G(F)$. Let $\cH(G)$ be the Hecke algebra of locally constant compactly supported complex valued functions on $G(F)$, where the product on this algebra is the one defined by the convolution integral with respect to the Haar measure giving the group $K$ measure $1$. We write $\cH_0(G)$ for the spherical Hecke algebra of $G$ with respect to $K$. We write $\rho$ for the half sum of the positive roots of $G$. 

We write $Z \subset G$ for the center of $G$, and we write $A \subset Z$ for the split center. Similarly $Z_M$ (resp. $Z_P$) is the center of the Levi-subgroup $M$ (resp. parabolic subgroup $P$); and we write $A_M$ (resp. $A_P$) for the split center of $M$. We write $A_0$ for $A_{P_0} \subset T$. We write $\ia_0 := X_*(A_0) \otimes \R$, and $C_0$ for the set of $x \in \ia_0$ such that for all roots $\alpha$ in $\Delta(A_0, \Lie(N_0))$ we have $\langle x, \alpha \rangle \geq 0$. 

Let $B(G)$ be the set of $\sigma$-conjugacy classes in $G(L)$, where $L$ is the completion of the maximal unramified extension of $F$ and $\sigma$ is the arithmetic Frobenius of $L$ over $F$. Let $\mu \in X_*(T)$ be a $G$-dominant minuscule cocharacter. Recall that Kottwitz has defined the subset $B(G, \mu) \subset B(G)$ of $\mu$-admissible isocrystals \cite{MR1485921, MR1411570}. 

Let $\D$ be the protorus over $F$ with character group given by $X_*(\D) = \Q$ and trivial Galois action. For any $b \in G(L)$ we have an unique morphism $\nu_b \colon \D_L \to G_L$ characterized by the following property: For every algebraic representation $(\rho, V)$ of $G$ on a finite dimensional vector space $V$ the composition $\rho \circ \nu_b$ determines the slope filtration on $(V \otimes L, \rho(b)(1 \otimes \sigma_L))$ \cite[\S 4]{MR809866}. Replacing $b$ by a $\sigma$-conjugate amounts to conjugating $\nu_b$ with some $G(L)$-conjugate. Moreover, one can replace $b$ so that $\nu_b$ has image inside the torus $A_{0, L}$, so that $\nu_b$ defines an element of $\ia_0$ \cite[p.~267]{MR1485921}\cite[1.7]{MR1411570}. Write $\li \nu_b$ for the unique element of $C_0$ whose orbit under the Weyl group meets $\nu_b$. The morphism $\li \nu_b$ is called the \emph{slope morphism} and the mapping $B(G) \to C_0, b \mapsto \li \nu_b$ is called the \emph{Newton map}. Note that the mapping $b \mapsto \li \nu_b$ is not injective in general (it is injective in case $G = \Gl_{n}(F)$).

Recall that we fixed an embedding $F \subset \lqp$. For each finite subextension $F' \subset \lqp$ of $F$ we have the unique mapping $H_T \colon T(F') \to X_*(T)_\R$ such that $q_{F'}^{-\langle \chi, H(t) \rangle} = |\chi(t)|$ for all $t \in T(F')$, where $q_{F'}$ is the cardinal of the residue field of $F'$, and the norm is normalized so that $|p|$ equals $q_{F'}^{-e}$ where $e$ is the ramification index of $F'/F$. By taking the union over all $F'$ we get a mapping $H_T \colon T(\lqp) \to X_*(T)_\R$. Consider the composition $H_A$ defined by $T(\lqp) \to X_*(T)_\R \to X_*(A)_\R = \ia_0$. Let $G(\lqp)_{\textup{ss}} \subset G(\lqp)$ be the subset of semisimple elements. If $g \in G(\lqp)_\textup{ss}$, then we may conjugate $g$ to an element $g'$ of $ T(\lqp)$ and then consider $H_A(g') \in \ia_0$. This element of $\ia_0$ is only defined up to conjugacy, but we can take a representative in the closed positive Weyl chamber $H(g) \in C_0^+$ which is well-defined. Thus we have a map $\Phi \colon G(\lqp)_{\textup{ss}} \to C_0$ defined on the semisimple elements. We extend the definition of $\Phi$ to $G(\lqp)$ by defining $\Phi(g) := \Phi(g_\textup{ss})$, where $g_\textup{ss}$ is the semisimple part of the element $g \in G(\lqp)$. We restrict to $ G(F) \subset G(\lqp)$ to obtain the mapping $\Phi \colon G \to C_0$. In Proposition~\ref{st:slopeproposition} we establish a relation between the map $\Phi$ and the Newton polygon mapping of isocrystals. 

We recall the definition of the norm $\cN$ of (certain) $\sigma$-conjugacy classes (cf. \cite{MR1007299} \cite[p.~799]{MR683003}). To any element $\delta \in G(F_\alpha)$ we associate the element $N(\delta) := \delta \sigma(\delta) \cdots \sigma^{\alpha - 1}(\delta) \in G(F_\alpha)$. For any element $\delta \in G(F_\alpha)$, defined up to $\sigma$-conjugacy, with semi-simple norm $N(\delta)$ one proves (see [\textit{loc. cit.}]) that $N(\delta)$ actually comes from a conjugacy class $\cN(\delta)$ in the group $G(F)$. 

\begin{proposition}\label{st:slopeproposition}
Let $\alpha$ be a positive integer and let $\delta \in G(F_\alpha)$ be an element of semi-simple norm, defined up to $\sigma$-conjugacy. Let $\gamma \in G(F)$ be an element in the conjugacy class $\cN(\delta)$, and let $b$ be the isocrystal with additional $G$-structure defined by $\delta$. Then $\li \nu_b = \alpha \cdot \Phi(\gamma) \in C_0$.
\end{proposition}
\begin{proof}
We first prove the case where $G$ is the general linear group. If $G = \Gl_{n,F}$, then an isocrystal ``with additional $G$-structure'' is simply an isocrystal, \ie a pair $(V, \Phi)$ where $V$ is an $n$-dimensional $L$ vector space and $\Phi$ is a $\sigma$-linear bijection from $V$ onto $V$. Because $b$ is induced by some $\delta \in G(F_\alpha)$, we find a $F_\alpha$-vector space $V'$ together with a $\sigma$-linear bijection $\Phi' \colon V' \to V'$ such that $(V, \Phi)$ is obtained from $(V, \Phi)$ by extending the scalars $V = V' \otimes_{F_\alpha} L$ and $\Phi(v' \otimes l) := \Phi'(v') \otimes \sigma(l)$. Then $(V', \Phi')$ is an $F_\alpha$-space in the terminology of Demazure \cite{MR883960}, and a theorem of Manin gives the relation $\li \nu_b = \alpha \cdot \Phi(\gamma)$ (cf. \cite[p.~90]{MR883960}). 

Drop the assumption that $G = \Gl_n$. Pick a representation $\rho \colon G \to \GL_V$ of $G$ in some finite dimensional $\qp$-vector space $V$. Then, by the statement for $\Gl_n$, we see that $\alpha \cdot \Phi_{\Gl_n}(\rho(\gamma))$ determines the slope filtration on the space $(V \otimes L, \rho(b)(1 \otimes \sigma_L))$. Thus $\rho \circ \li \nu_b = \alpha \cdot \Phi_{\Gl_n}(\rho(\gamma))$ for all $\rho$, and then the equality is also true for the group $G$. 
\end{proof}

\subsection{Truncated traces}
In this section we introduce the concept of truncated traces of smooth representations with respect to elements of the set $B(G)$, \ie the isocrystals with additional $G$-structure. We will then compute these truncated traces on the Steinberg representation and on the trivial representation. 

Using the mapping $\Phi$ from the previous section we define the truncated traces with respect to an arbitrary element $b \in B(G)$:

\begin{definition}
Let $\nu \in C_0$. We write $\Omega_\nu^G$ for the set of $g \in G$ such that there exists a $\lambda \in \R_{>0}$ such that $\Phi(g) = \lambda \cdot \nu \in C_0$. We let $\chi_\nu^G$ be the characteristic function on of the subset $\Omega_\nu^G$ of $G$. Let $b \in B(G)$ be an isocrystal with additional $G$-structure. Then we will write $\chi_b^G := \chi_{\li \nu_b}^G$ and $\Omega_b^G := \Omega_{\li \nu_b}^G$.
\end{definition}

\begin{remark}
The Newton mapping $B(G) \owns b \mapsto \li \nu_b \in C_0$ is \emph{injective} for a simply connected, connected quasi-split reductive group over a non-Archimedean local field \cite[\S 6]{MR1485921}. 
\end{remark}

Let $P = MN$ be a standard parabolic subgroup of $G$ and let $A_P$ be the split center of $P$, we write $\eps_P = (-1)^{\dim(A_P/A_G)}$. To the parabolic subgroup $P$ we associate the subset $\Delta_P \subset \Delta$ consisting of those roots acting non trivially on $A_P$. Define $\ia_P$ to be $X_*(A_P)_\R$, define $\ia_P^G$ to be the quotient of $\ia_P$ by $\ia_G$, and define $\ia_P^+$ to be the set of $x \in \ia_P$ such that for all roots $\alpha$ in $\Delta_P$ we have $\langle x, \alpha \rangle > 0$.

We recall the definition of the obtuse and acute Weyl-chambers \cite{MR1989693, waldlabessetwisted}. Let $P$ be a standard parabolic subgroup of $G$. We write $\ia_0 = \ia_{P_0}$ and $\ia_0^G = \ia_{P_0}^G$. For each root $\alpha$ in $\Delta$ we have a coroot $\alpha^\vee$ in $\ia_0^G$. For $\alpha \in \Delta_P \subset \Delta$ we send the coroot $\alpha^\vee \in \ia_0^G$ to the space $\ia_P^G$ via the canonical surjection $\ia_0^G \surjects \ia_P^G$. The set of these restricted coroots $\alpha^\vee|_{\ia_P^G}$ with $\alpha$ ranging over $\Delta_P$ form a basis of the vector space $\ia_P^G$. By definition the set of fundamental weights $\{ \varpi_\alpha \in \ia_P^{G*}\,\,|\,\, \alpha \in \Delta_P\}$ is the basis of $\ia_P^{G*} = \Hom(\ia_P^G , \R)$ dual to the basis $\{\alpha^\vee_{\ia_P^G}\}$ of coroots. The \emph{acute Weyl chamber} $\ia_P^{G+}$ is the set of $x \in \ia_P^G$ such that for all $\alpha \in \Delta_P$ we have $\langle \alpha,x \rangle > 0$. The \emph{obtuse Weyl chamber} ${}^+\ia_P^G$ is the set of $x \in \ia_P^G$ such that for all $\alpha \in \Delta_P$ we have $\langle \varpi_\alpha^G, x \rangle > 0$. We let $\tau_P^G$ (resp. $\widehat \tau_P^G$) be the characteristic function on the space $\ia_P^G$ of the acute (resp. obtuse) Weyl chamber. We define the function $\chi_N$ to be the composition $\tau_P^G \circ (\ia_P \surjects \ia_P^G) \circ H_M$, and we define the function $\widehat \chi_N$ to be the composition $\widehat \tau_P^G \circ (\ia_P \surjects \ia_P^G) \circ H_M$. The functions $\chi_N$ and $\widehat \chi_N$ are locally constant and $K_M$-invariant, where $K_M = M(\cO_F)$.

Let $b \in B(G)$ be an isocrystal with additional $G$ structure and let $\li \nu_b \in C_0$ be its slope morphism. For any standard parabolic subgroup $P \subset G$ we have the subset $\ia_P^+ \subset C_0$. Let $P_b$ be the standard parabolic subgroup of $G$ such that $\li \nu_b \in \ia_{P_b}^+$. We call the group $P_b$ the subgroup of $G$ \emph{contracted} by the isocrystal $b \in B(G)$. These groups are precisely the parabolic subgroups appearing in the Kottwitz decomposition of the set $B(G)$ (see \cite[5.1.1]{MR1485921}). We write $P_b = M_b N_b$ for the standard decomposition of $P_b$.

We write $\pi_{0 P}$ for the projection from the space $\ia_0$ onto $\ia_P$, it sends an element $X \in \ia_0$ to its average under the action of the Weyl group.

We introduce a certain characteristic function on $G$ associated to the isocrystal $b \in B(G)$:

\begin{definition}
Let $P_b = M_b N_b$ be the standard parabolic subgroup of $G$ contracted by $b$. We define $\eta_b$ to be the characteristic function on $G$ of the set of elements $g \in G$ such that there exists a $\lambda \in \R^\times_{>0}$ such that $\pi_{0P} (\Phi(g)) = \lambda \cdot \li \nu_b \in \ia_P^+$. 
\end{definition}

\begin{remark}
If the isocrystal $b$ is basic, then we have $P = G$, and the element $\li \nu_b \in C_0$ is central. Therefore the function $\eta_b$ is \emph{spherical}. 
\end{remark}

In case the isocrystal $b \in b(G)$ is basic then $\chi_b^G$ coincides with $\eta_b\chi_c^G$:

\begin{lemma}\label{st:basicstruncatedcompact}
Let $b \in B(G)$ be a basic isocrystal. Then we have $\chi_b^G = \eta_b \chi_c^G$. 
\end{lemma}
\begin{proof}
Let $g \in G$, and consider $\Phi(g) \in C_0$. Then $g$ is compact if and only if it contracts $G$ as parabolic subgroup (meaning that $\Phi(g)$ lies in $\ia_G \subset C_0^+$). Assume $g$ is compact. Then $\chi_b^G(g) = 1$ if and only if the slope morphism $\li \nu_b$ of $b$ lies in $\ia_G$, \ie if and only if the centralizer of the slope morphism of $b$ is equal to $G$. But that means that $b$ is basic. Conversely, assume $b$ is basic. Then its slope morphism is central, thus $\chi_c^G(g) = 1$ if and only if $g$ contracts $G$, \ie $g$ is compact. Furthermore we have $\eta_b(g) = 1$ because $\Phi(g)$ equals $\li \nu_b$ up to a positive scalar. This completes the proof. 
\end{proof}

We call the collection of subsets $\Omega_b^G$ for $b \in B(G)$ the \emph{Newton polygon stratification} of the group $G$. For our proofs we will also need to study another stratification, called the \emph{Casselman stratification} of $G$:

\begin{definition}
Let $Q$ be a standard parabolic subgroup of $G$. We define $\Omega_Q^G \subset G$ to be the subset of elements $g \in G$ contracting \cite[\S 1]{MR1068388} a parabolic subgroup conjugate to $Q$. Write $\chi_Q^G$ for the characteristic function on $G$ of the subset $\Omega_Q^G \subset G$. These sets $\Omega_Q^G$ form the Casselman stratification of $G$.
\end{definition}

For truncated traces with respect to the Casselman stratification we have:

\begin{proposition}\label{st:casselnewton}
Let $Q = LU$ be a standard parabolic subgroup of $G$. Let $f \in \cH(G)$ be a locally constant function with compact support. Then we have $\Tr(\chi_Q^G f, \pi) = \Tr(\chi_U \chi_c^L \li f^{(Q)}, \pi_U(\delta_{Q}^{-1/2}))$. 
\end{proposition}
\begin{proof}
By the Proposition \cite[prop 1.1]{MR1068388} on compact traces, for all functions $f$ on $G$, the full trace $\Tr(f, \pi)$ is equal to the sum of compact traces $\sum \Tr_{M}(\chi_c^M \li f^{(P)}, \pi_N(\delta_P^{-1/2}))$, where the sum ranges over the standard parabolic subgroups $P = MN$ of $G$. Consider only those functions of the form $\chi_Q^G f \in \cH(G)$. Then we obtain that the trace $\Tr(\chi_Q^G f, \pi)$ is equal to the sum $\sum \Tr_{M}(\chi_c^M \chi_Q^G \li f^{(P)}, \pi_N(\delta_P^{-1/2}))$ where $P = MN$ ranges over the standard parabolic subgroups of $G$. Observe that $\chi_c^M \chi_Q^G = 0$ if $P \neq Q$. Therefore only the term corresponding to $P = Q$ remains in the sum. This completes the proof.
\end{proof}

Let us now explain the relation between the Casselman stratification and the Newton stratification. The following Proposition gives the relation between the Casselman stratification of $G$ and the Newton stratification:

\begin{proposition}
For all $b \in B(G)$ we have $\Omega_b^G \subset \Omega_{P_b}^G$. 
\end{proposition}
\begin{proof}
Assume that $g \in \Omega_b^G$. Then $\Phi(g) = \lambda \li \nu_b \in \ia_0$. Let $P$ be the standard parabolic subgroup of $G$ conjugate to the parabolic subgroup of $G$ contracted by $g$. Then $\li \nu_b = \lambda \Phi(g) \in \ia_P^+$. Then, by definition, $P$ is the parabolic subgroup contracted by $b$. This completes the proof.
\end{proof}

\begin{example}
The inclusion $\Omega_b^G \subset \Omega_{P_b}^G$ is strict in general. Consider for example the case $G = \Gl_{n,\qp}$ to see that it is non-strict only in particular cases, such as when $n = 2$. In the particular case of the Shimura varieties of Harris-Taylor \cite{MR1876802}, the Casselman stratification also separates the isocrystals. 
\end{example}

We will now compute the truncated trace on the Steinberg representation. 

\begin{definition}
Let $\xi_b^{St}$ be the characteristic function on $T$ defined by $\xi_b^{St} := \widehat \chi_{N_0 \cap M_b} \chi_{N_b} \eta_{b}$, where with the notation $\widehat \chi_{N_0 \cap M_b}$ we mean the characteristic function on the Levi subgroup $M_b \subset G$, corresponding to the obtuse chamber relative to the minimal parabolic subgroup of $M_b$.
\end{definition}

\begin{proposition}\label{st:truncatedtraceonsteinberg} Let $f \in \cH_0(G)$ be a spherical Hecke operator. Then we have
\begin{equation}\label{eqn:Trace_Steinberg_Abstract}
\Tr(\chi_b^G f, \St_G) = \eps_{P_0 \cap M_b} \Tr_T(\xi_b^{St} f^{(P_0)}, \one(\delta_{P_b}^{-1/2} \delta_{P_0 \cap M_b}^{1/2})).
\end{equation}
\end{proposition}
\begin{proof} Write $P = MN$ for the parabolic subgroup contracted by the isocrystal $b$. We compute: 
\begin{align}\label{eqn:stpropA}
\Tr(\chi_b^G f, \St_G) = \Tr_M(\chi_b^G \chi_N f^{(P)}, (\St_G)_N(\delta_P^{-1/2})),
\end{align}
(Proposition~\ref{st:casselnewton}). Let $b_M \in B(M)$ be a $G$-regular basic element such that its image in $B(G)$ is equal to $b$ \cite[prop. 6.3]{MR809866}. By [\textit{loc. cit.}] the set of all such $b_M$ are $G$-conjugate. As functions on $M$ we have $\chi_b^G \chi_N = \chi_{b_M}^M \chi_N$. Therefore we may simplify Equation~\eqref{eqn:stpropA} to $\Tr_M(\chi_{b_M}^M \chi_N f^{(P)}, (\St_G)_N(\delta_P^{-1/2}))$. By Lemma~\ref{st:basicstruncatedcompact} the latter trace equals $\Tr_M(\chi_{c}^M \eta_b \chi_N f^{(P)}, (\St_G)_N(\delta_P^{-1/2}))$. In our article \cite{kret1} we computed the compact traces on the Steinberg representation for all spherical Hecke operators. By \cite[Prop.~1.13]{kret1} we get
\begin{align}\label{eqn:replacethis}
\Tr(\chi_b^G & f, \St_G) = \eps_{P_0 \cap M} \Tr(\widehat \chi_{N_0 \cap M} \eta_b \chi_N f^{(P_0)}, \one(\delta_P^{-1/2} \delta_{P_0 \cap M}^{1/2})).
\end{align} 
This completes the proof. 
\end{proof}

In the same way one may compute the truncated traces on the trivial representation. We have to introduce two more notations. Let $\widehat \chi^{\leq}_{N_0 \cap M_b}$ be the characteristic function on $M_b$ corresponding to the negative closed obtuse chamber in $\ia_P$. Then we define:

\begin{definition}
Let $b \in B(G)$ be an isocrystal. We define $\xi_b^{\one} := \widehat \chi^{\leq}_{N_0 \cap M_b} \chi_{N_b} \eta_b$. 
\end{definition} 

\begin{proposition}\label{st:trivialprop} We have
$\Tr(\chi_b^G f, \one) = \Tr_T(\xi_b^{\one} f^{(P_0)}, \one(\delta_{P_0}^{-1/2}))$.
\end{proposition}
\begin{proof}
The proof of Proposition~\ref{st:truncatedtraceonsteinberg} may be repeated without change up to Equation~\eqref{eqn:replacethis}. Replace the result in that last Equation with the result from Proposition 3.1 from \cite{kret2}. This Proposition gives the compact trace on the trivial representation for any Hecke operator (and any unramified group).
\end{proof}

\begin{remark}
With a method similar to the above one may compute the truncated traces on the irreducible subquotients of the $G$-representation on the space $C^\infty(G/P_0)$ of locally constant functions on $G/P_0$. 
\end{remark}

\subsection{The function of Kottwitz}
Let $G$ be a connected, reductive unramified group over $\qp$, let $P_0$ be a Borel subgroup of $G$. Let $T$ be the Levi-component of $P_0$. Then $T$ is a maximal torus in $G$, and let $W$ be the absolute Weyl group of $T$ in $G$. Let $\mu \in X_*(T)$ be a minuscule cocharacter. 

We write in this section $E$ for an arbitrary, finite unramified extension of $\qp$. In later sections, the field $E$ that we consider here will be the completion of the reflex field at a prime of good reduction. We fix an embedding of $E$ into $\lqp$, and for each positive integer $\alpha$ we write $E_\alpha \subset \lqp$ for the unramified extension of degree $\alpha$ of $E$.

\begin{definition}\label{def:kottwitzfunction} (cf. \cite{MR761308}). Let $\alpha$ be a positive integer, and $E_\alpha$ the unramified extension of $E$ of degree $\alpha$ contained in $\lqp$. We write $W_\alpha$ for the subgroup $W(G(E_\alpha), T(E_\alpha))$ of $W$. Write $S_\alpha$ for a maximal $E_\alpha$-split subtorus of $G_{E_\alpha}$. We define $\phi_{G, \mu, \alpha} \in \cH_0(G(E_\alpha))$ to be the spherical function whose Satake transform is equal to 
\begin{equation}\label{eqn:kottwitzfunction}
p^{-\alpha \langle \rho_G, \mu \rangle} \sum_{w \in W_\alpha / \stab_{W_\alpha}(\mu)} [w(\mu)] \in \C[X_*(S_\alpha)]^{W_\alpha},
\end{equation}
where $\stab_{W_\alpha}(\mu) \subset W_\alpha$ is the stabilizer of $\mu$ in the group $W_\alpha$. We define $f_{G, \mu, \alpha}$ to be the function obtained from $\phi_{G, \mu, \alpha}$ via base change from the group $G(E_\alpha)$ to the group $G(\qp)$. We call $f_{G, \mu,\alpha}$ the \emph{function of Kottwitz}.
\end{definition}

\begin{remark}
Kottwitz proves in \cite{MR761308} that the definition of the Kottwitz functions $f_{G, \mu, \alpha}$ and $\phi_{\G, \mu, \alpha}$ coincides with the definition that we gave at the end of section 2.
\end{remark}

\begin{remark}
We note that the notation for the functions $f_{G, \mu, \alpha}$ and $\phi_{G, \mu, \alpha}$ is slightly abusive, as they also depend on the field $E$. Because confusion will not be possible, we have decided to drop the field $E$ from the notations.
\end{remark}

\begin{proposition}\label{st:constantterms}
Let $P = MN$ be a standard parabolic subgroup of $G$. We have
$$
f_{G,\mu, \alpha}^{(P)} = q^{-\alpha \langle \rho_G - \rho_M, \mu \rangle} \sum_{w \in W_\alpha/\stab_{W_\alpha}(\mu) W_{M, \alpha}} f_{M,w(\mu),\alpha} \in \cH_0(M),
$$
where $\stab_{W_\alpha}(\mu) W_{M, \alpha} \subset W_\alpha$ is the subgroup of $W_\alpha$ generated by the group $W_{M, \alpha}$ of the Weyl group of $T(E_\alpha)$ in $M(E_\alpha)$ and the stabilizer subgroup of $\mu$ in $W_\alpha$.
\end{proposition}
\begin{proof}
Compute the Satake transform of both sides to see that they are equal.
\end{proof}

The integer $\alpha$ is the degree of the finite field over which we count points in the Newton stratum. We only want to show that the Newton-strata are non-empty. Therefore, we will often take $\alpha$ large so that the combinatorial problems simplify (large in the divisible sense).


\section{Unitary groups}
We prove that the Newton strata of unitary Shimura varieties of PEL type are nonempty. We continue with the local notations from the previous section, in particular $F$ is a local field.

\subsection{Unitary isocrystals} Let $G$ be an unramified unitary group over $F$ splitting over the extension $F_2/F$. 
The absolute root system of $G$ is isomorphic to the usual root system in $\R^n$ of type A (cf. Bourbaki \cite[chap. 6]{MR1890629}), and the non-trivial element of the group $\Gal(F_2/F)$ acts on $\R^n$ via the operator $\theta$ defined by $(x_1, x_2, \ldots, x_n) \longmapsto (-x_n, -x_{n-1}, \ldots, -x_1)$. The space $\ia_0$ is the subspace of $\theta$ invariant elements in $\R^n$, it equals the set of $(x_i) \in \R^n$ with $x_i = -x_{n+1 -i}$ for all indices $i$. The dimension of $\ia_0$ equals $\lfloor n/2 \rfloor$.

Whenever $b \in B(G)$ is an isocrystal with $G$-structure, we have its slope morphism $\li \nu_b \in C_0$. We may view the slope morphism $\li \nu_b$ as an $\theta$-invariant element of $\R^n$. This way we get the slopes $\lambda_1, \lambda_2, \ldots, \lambda_n$ of $b$. These slopes are just the coordinates of the vector $\li \nu_b \in \R^n$. We order them so that $\lambda_1 \leq \lambda_2 \leq \cdots \leq \lambda_n$. The slopes satisfy the property $\lambda_i = -\lambda_{n + 1 - i}$. 

We associate to the slopes $\lambda_i$ the Newton polygon $\cG_b$ of $b$. The Newton polygon is by definition the continuous piecewise linear function from the real interval $[0, n]$ to $\R$ with the property that the only points where it is possibly not differentiable are the integral points $[0, n] \cap \Z$; the value of $\cG_b$ at these points is defined by: $\cG_b(0) := 0$ and $\cG_b(i) := \lambda_1 + \lambda_2 + \ldots + \lambda_i$. 

Due to the $\theta$-invariance, we have $\cG_b(n) = \lambda_1 + \ldots + \lambda_n = 0$. 
Furthermore the graph (or polygon) $\cG_b$ is symmetric around the vertical line that goes through the point $(\frac n2, 0)$. In Figure 1 we show a typical unitary Newton polygon. In particular negative slopes may occur, which does not happen for the group $\Gl_n(F)$ nor for the group $\Gsp_{2g}(F)$.

\begin{figure}
\begin{center}
\input{Figure_unitary.tex}
\end{center}
\begin{caption}{The dark line is an example of the Newton polygon of an isocrystal $b$ with additional $U_{10}^*$-structure. The horizontal line from $(0,0)$ to $(10, 0)$ is the Newton polygon of the basic isocrystal. The vertical dotted line indicates the mirror symmetry of the Newton polygons of the $G$-isocrystals.}
\end{caption}
\end{figure}

Let us now determine what the Hodge polygons looks like. The minuscule cocharacter $\mu$ is defined over $\li F$, and is given by
$$
\mu_s = (\underset{n-s}{\underbrace{0, 0, \ldots, 0}}, \underset{s}{\underbrace{1, 1, \ldots, 1}}) \in \Z^n \subset \R^n, 
$$ 
for some integer $s$ with $0 \leq s \leq n$. To define the set $B(G, \mu)$ Kottwitz \cite[\S 6]{MR1485921} takes the average of $\mu$ under the Galois action to get 
$$
\li \mu := \tfrac 12 (\mu + \theta(\mu)) = (\underset{s'}{\underbrace{-\tfrac 12, -\tfrac 12, \ldots, -\tfrac 12}}, \underset{n-2s'}{\underbrace{0, 0, \ldots, 0}}, \underset{s'}{\underbrace{\tfrac 12, \tfrac 12, \ldots, \tfrac 12}}) \in \ia_0 \subset \R^n,
$$
where $s' := \min(s, n-s)$. To this element $\li \mu \in \R^n$ we may associate in the same manner a graph $\cG_\mu$ as in Figure $1$. Then $b \in B(G)$ lies in $B(G, \mu)$ if and only if the end point of $\cG_b$ is $(n,0)$ and if $\cG_b$ lies above\footnote{Lies above in the non-strict sense, the two graphs may touch, or even be the same (the ordinary case).} the graph $\cG_{\mu}$. 

\subsection{PEL datum}\label{sect:PEL_datum}
Let $G/\Q$ be a unitary group of similitudes arising from a PEL type Shimura datum \cite[\S 5]{MR1124982}. We recall briefly the definition of $G$ from [\textit{loc. cit.}]. Consider the following notations:
\begin{itemize}
\item $B/\Q$ is a finite dimensional simple algebra;
\item $F$ is the center of $B$, by assumption $F$ is a CM field; 
\item $*$ is a positive involution on $B$ over $\Q$ inducing complex conjugation on $F$; 
\item $F^+ \subset F$ is the totally real subfield of $F$; 
\item $V$ is a nonzero finitely generated left $B$-module;
\item $(\cdot, \cdot)$ is a non-degenerate $\Q$-valued alternating form on $V$ such that $(bv, w) = (v, b^* w)$ for all $v, w \in V$ and all $b \in B$. 
\item $h \colon \C \to B_\R^{\textup{opp}}$ a morphism of real algebras, such that $h(\li z) = z^*$ for all complex numbers $z$ and the involution $x \mapsto h(i)x^*h^(i)^{-1}$ is positive;
\item $G/\Q$ is the algebraic group such that for all commutative $\Q$-algebras $R$ the set $G(R)$ is equal to the set of $g \in \End_B(V)^\times$ such that there exists $c(g) \in R^\times$ such that $(g \cdot, g \cdot) = c(g) (\cdot, \cdot)$ on $V$; 
\item finally $X$ is the $G(\R)$-conjugacy class of the morphism $h$.
\end{itemize}
The couple $(G, X)$ is a Shimura-datum of PEL-type. We write furthermore:
\begin{itemize}
\item $G_1 \subset G$ for the kernel of the similitudes ratio $G \surjects \Gm$; 
\item $G_1$ is obtained by restriction of scalars of a unitary group defined over the totally real field $F^+$ and  we write $G_0$ for this group. 
\end{itemize}

The group $G_{1,\qp}$ is isomorphic to the product $G_{1, \qp} \cong \prod_{\wp | p} G_{1,\wp}$
where $\wp$ ranges over the $F^+$-places above $p$. For each $\wp$ the group $G_{1,\wp}$ is either the restriction of scalars to $\qp$ of $\Gl_{n, F_{\wp}^+}$ or of an unramified unitary group over $F_{\wp}^+$. 

We will often study the group $G_{1,\qp}$ factor by factor. Thus, in this article we do not  owrk only with unramified unitary groups, but with the slightly more general class of groups of the form $\Res_{F'_\wp/F_\wp} U$, where $F'_\wp/F_\wp$ is some unramified extension and $U$ is an unramified unitary group over $F_\wp'$. The study of isocrystals over these groups reduces quickly to the study of isocrystals over the group $U$ (which we did above), by the Shapiro bijection (cf. \cite[6.5.3]{MR1485921}): $B(\Res_{F'_\wp/F_\wp} U) = B_{F'_{\wp}}(U)$,
(the subscript ``$F'_{\wp}$'' indicates that we work with $\sigma'$-conjugacy classes, where $\sigma'$ is the arithmetic Frobenius of $\lqp$ over $F_{\wp}'$). Under the Shapiro bijection the subset $B(\Res_{F'_\wp/F_\wp} U, \mu_\wp)$ corresponds to the subset $B_{F_\wp'}(U, \mu'_\wp)$ of $B_{F'_{\wp}}(U)$, where $\mu'_\wp = \sum_{v \in V(\wp)} \mu_{s_v}$. Thus, combinatorics for isocrystals with $\Res_{F'/F} U$-structure is almost the same as the combinatorics for the case $F' = F$; only the Hodge polygons are slightly more complicated.

We recall how the functions of Kottwitz $\phi_\alpha$ and $f_\alpha$ are constructed~\cite[\S 5]{MR1124982} \cite[p.~173]{MR1044820}. Let $E$ be the reflex field, let $\p$ be an $E$-prime where the Shimura variety has good reduction in the sense of \cite[\S 6]{MR1124982}. In particular, we have the following notations/facts:
\begin{itemize} 
\item The field $E$ is unramified at $\p$; 
\item $p$ is the rational prime number lying below $\p$;
\item $E_\p$ is the completion of $E$ at $\p$; 
\item fix once and for all an embedding $E_\p \subset \lqp$; 
\item let for each positive integer $\alpha$, the field $E_{\p, \alpha} \subset \lqp$ be the unramified extension of $E_\p$ of degree $\alpha$. 
\end{itemize}
In the PEL datum there is fixed a $*$-morphism $h \colon \C \to \End(B)_\R^{\textup{opp}}$. This morphism induces a morphism of algebraic groups from Deligne's torus $\Res_{\C/\R} {\mathbb G}_{\textup{m}}$ to the group $G_\R$. Tensor this morphism with $\C$ to get a morphism from ${\mathbb G}_{\textup{m}} \times {\mathbb G}_{\textup{m}}$ to $G_\C$ and then restrict to the factor ${\mathbb G}_{\textup{m}}$ of the product ${\mathbb G}_{\textup{m}} \times {\mathbb G}_{\textup{m}}$ corresponding to the identity $\R$-isomorphism $\C \to \C$. This way we obtain a cocharacter $\mu \in X_*(G)$. We quote from Kottwitz's article at Ann Arbor, p.~173: The $G(\C)$ conjugacy class of $\mu$ gives a $G(\lqp)$ conjugacy class of morphisms fixed by the Galois group $\Gal(\lqp/E_{\p, \alpha})$. Let $S_\alpha$ be a maximal $E_{\p, \alpha}$-split torus in the group $G$ over the ring of integers $\cO_{E_{\p, \alpha}}$. Using Lemma (1.1.3) of \cite{MR761308} we choose $\mu$ so that it factors through $S_\alpha$. Then $\phi_\alpha = \phi_{G, \mu, \alpha}$ is the characteristic function of the double coset $G(\cO_{\p, \alpha}) \mu(p^{-1})G(\cO_{\p, \alpha})$. The function $f_\alpha = f_{G, \mu, \alpha}$ is by definition the base change \cite{MR1007299, MR2856377} of $\phi_\alpha$ from the group $G(E_{\p, \alpha})$ to the group $G(\qp)$.

\subsection{The class of $\iR(b)$-representations}
For the global applications to Shimura varieties we define a class representations $\iR_1(b)$ of positive Plancherel density on which the truncated trace of the Kottwitz functions are non-zero. In fact we take for \emph{most} of the isocrystals $b \in B(G, \mu)$ simply the Steinberg representation at $p$, but there are some exceptions where the truncated trace on the Steinberg representation vanishes; in those cases we take a different representation.

We make the function of Kottwitz explicit in case $G$ is either the restriction of scalars of a general linear group over $F^+$ or the restriction of scalars of an unramified unitary group over $F^+$. 

We assume that we are in one of the following two cases:
\begin{equation}\label{eqn:Local_group_Unitary}
G = \begin{cases} \Res_{F^+/\qp}(\Gl_{n, F^+}) & \textup{(linear type)} \cr
\Res_{F^+/\qp}(U) & \textup{(unitary type)} 
\end{cases}
\end{equation}
where $F^+/\qp$ is a finite unramified extension, and where $U/F^+$ is an unramified unitary group, outer form of $\Gl_{n, F^+}$. These groups $G$ occur as the components in the product decomposition of $G_\qp$. We assume that the cocharacter $\mu \in X_*(T)$ arises from a PEL-type datum. 

We begin with the linear case. We have a cocharacter $\mu \in X_*(T)$. Thus, for each $\qp$-embedding $v$ of $F^+$ into $\lqp$ we get a cocharacter $\mu_v$ of the form
$$
(\underset{s_v}{\underbrace{1, 1, \ldots, 1}}, \underset{n-s_v}{\underbrace{0, 0, \ldots, 0}}) \in \Z^n.
$$ 
To each such integer $s_v$ we associate the spherical function $f_{n\alpha s_v}$ on $\Gl_n(F^+)$ whose Satake transform is defined by
\begin{equation}\label{eqn:Simple_Kottwitz_function_GL_f_nas}
\cS(f_{n\alpha s_v}) = q^{\frac {s(n-s)}2 \alpha} \sum_{i_1 < i_2 < \ldots < i_{s_v}} X_{i_1}^\alpha X_{i_2}^\alpha \cdots X_{i_{s_v}}^\alpha \in \C[X_1^{\pm 1}, \ldots, X_n^{\pm 1}]. 
\end{equation}
We write $V_\alpha$ for the set of $\Gal(\lqp/E_{\alpha})$-orbits in the set $\Hom(F^+, \lqp)$. If $v \in V_\alpha$ is such an orbit, then this orbit corresponds to a certain finite unramified extension $E_{\alpha}[v]$ of $E_{\alpha}$. Let $\alpha_v$ be the degree over $\qp$ of the field $E_{\alpha}[v]$, we then have $E_{\alpha}[v] = E_{\alpha_v}$. The function $f_\alpha$ is given by
\begin{equation}\label{eqn:Kottwitz_function_GLN}
f_\alpha = \prod_{v \in V_\alpha} f^{\Gl_n(F^+)}_{n \alpha_v s_v} \in \cH_0(G(\qp)),
\end{equation}
where the product is the convolution product (cf. \cite[Prop. 3.3]{kret1}). 

Let us now assume that we are in the unitary case (cf. Equation~\eqref{eqn:Local_group_Unitary}). We will make the function $f_{G, \mu, \alpha}$ explicit only in case $\alpha$ is even. To obtain the function of Kottwitz on $G$, we have to apply base change from $G(E_\alpha)$ to $G(\qp)$. Assume that $\alpha$ is even. Let $\Q_{p^2}$ be the quadratic unramified extension of $\qp$ contained in $\lqp$. The base change factors over the composition of base changes $G(E_\alpha) \rightsquigarrow G(\Q_{p^2}) \rightsquigarrow G(\qp)$. The base change of $\phi_\alpha$ to $G(\Q_{p^2})$ is a function of the form $f_{G^+, \mu, (\alpha/2)}$ on the group $G^+ = \Res_{\Q_{p^2}/\qp}(G_{\Q_{p^2}})$. Explicitly, the last quadratic base change $G(\Q_{p^2}) \rightsquigarrow G(\qp)$ is given by:
\begin{align}\label{eqn:Base_Change_Map_Unitary_Group}
\Psi \colon \C[X_1^{\pm 1}, \ldots, X_n^{\pm 1}]^{\iS_n} &\to \C[X_1^{\pm 1}, \ldots, X_n^{\pm 1}]^{\iS_m \rtimes (\Z/2\Z)^m}, \cr
 X_i & \longmapsto \begin{cases} X_i & 1 \leq i \leq \lfloor n/2 \rfloor, \cr
1 & i = \lfloor \tfrac n2 + 1 \rfloor, \textup{ and $n$ is odd,} \cr
X_{n + 1 -i}^{-1} & n+1 - \lfloor n/2 \rfloor \leq i \leq n,
\end{cases}
\end{align}
where $m := \lfloor \tfrac n2 \rfloor$ (cf. \cite{MR2856377}). Thus we get $f_{G, \mu, \alpha} = \Psi f_{G^+, \mu, (\alpha/2)}$.

\begin{lemma}\label{st:prenonzerolemma} Let $G$ be an algebraic group over $\qp$ defined as in Equation~\eqref{eqn:Local_group_Unitary}. Let $\pi$ be a generic unramified representation of $G$, and $f = f_{G,\mu,\alpha}$ a function of Kottwitz, and $b \in B(G)$ an isocrystal. Let $\alpha \in \Z_{>0}$ be an integer, sufficiently divisible such that $W_\alpha$ is the absolute the Weyl group of $T$ in $G$. Then, the truncated trace $\Tr(\chi_b^G f_{G,\mu,\alpha}, \pi)$ is non-zero if and only if there exists some $w \in W_\alpha$ and some $\lambda \in \R^\times_{>0}$ such that $w(\mu) = \lambda \li \nu_b \in \ia_0^G$.
\end{lemma}
\begin{remark}
In case $G$ is the general linear group, then there exists a pair $w \in W, \lambda \in \R^\times_{>0}$ such that $w(\mu) = \lambda \li \nu_b$ if and only if the slopes $\lambda_i$ of $b$ all lie in the set $\{0, 1\}$. 
\end{remark}
\begin{proof}
We have $\pi = \Ind_T^G(\rho)$, where $\rho$ is some smooth character of the torus $T$. By van Dijk's formula for truncated traces \cite[Prop. 1.1]{kret1} we have $\Tr(\chi_b^G f, \pi) = \Tr(\chi_b^G f^{(P_0)}, \rho)$. The truncation operation $h \mapsto \chi_b^G h$ on $\cH_0(T)$, corresponds via the Satake transform to an operation on $\C[X_*(T)]$ sending certain monomials $[M] \in \C[X_*(T)]$ associated to elements $M \in X_*(T)$ to zero, and leaves certain other monomials invariant. Thus to compute the trace $\Tr(\chi_b^G f^{(P_0)}, \rho)$ one takes the set of monomials $[w(\mu)]$, $w \in W$ occurring in $f^{(P_0)}$, and removes some of them (maybe all), and then evaluate those which are left at the Hecke matrix of $\rho$. The lemma now follows from the observation that $\chi_b^G \cS^{-1}_T [w(\mu)] \neq 0$ if and only if $w(\mu) = \lambda \li \nu_b$ for some positive scalar $\lambda \in \R_{>0}^\times$. This completes the proof.
\end{proof}

We have to distinguish further between (essentially) two cases at $p$. The case the group is the general linear group, and the case where the group is the unramified unitary group. We begin with the general linear group.

\begin{proposition}\label{st:ACclassicalgroups}
Let $G$ be an algebraic group over $\qp$ defined as in Equation~\eqref{eqn:Local_group_Unitary}, and assume it is of linear type, so $G(\qp) = \Gl_n(F^+)$. Let $b \in B(G, \mu)$ be a $\mu$-admissible isocrystal having the property that the number of slopes equal to $0$ is at most $1$, and the number of slopes equal to $1$ is also at most $1$. Let $\chi$ be an unramified character of $\Gl_n(F^+)$. Then, for $\alpha$ sufficiently divisible, we have $\Tr(\chi_b^G f_{G,\mu,\alpha}, \St_G(\chi)) \neq 0$.
\end{proposition}
\begin{remark}
In the proof of the Proposition we use the divisibility of $\alpha$ at two places. First, it simplifies the function of Kottwitz (cf. Equation~\eqref{eqn:Kottwitz_function_GLN}). Second, we want $\alpha$ sufficiently divisible so that the Weyl group $W(T(E_\alpha), G(E_\alpha))$ relative to the field $E_\alpha$ is the full Weyl group. 
\end{remark}
\begin{remark}
In case the isocrystal $b$ has two or more slopes with value $0$ (or $1$), then the truncated trace of the Kottwitz function on the Steinberg representation vanishes. 
\end{remark}
\begin{proof}
By Proposition~\ref{st:truncatedtraceonsteinberg} we have to show that the function $\xi_b^{\St} f_{G, \mu, \alpha}^{(P_0)}$ does not vanish. Recall that the function $f_{G, \mu, \alpha}$ is obtained from a function $\phi_\alpha$ through base change from the group $\Gl_n(F^+ \otimes E_{\alpha})$. Let us first assume that the $E_\alpha$-algebra $F^+ \otimes E_{\alpha}$ is a field. In that case we have that $f_{G, \mu, \alpha} = f_{n\alpha s}$ in the notations from \cite[p.~10]{kret1}, \ie $\cS(f_{G, \mu, \alpha})$ is (up to scalar) an elementary symmetric function in the Satake algebra, 
\begin{equation}\label{eqn:sattransform}
\cS(f_{G, \mu, \alpha}) = q^{\frac {s_v(n-s_v)}2 \alpha} \sum_{i_1 < i_2 < \ldots < i_s} X_{i_1}^{d \alpha} X_{i_2}^{d \alpha} \cdots X_{i_s}^\alpha \in \C[X_1^{\pm 1}, \ldots, X_n^{\pm 1}]. 
\end{equation}
We have to show that under the truncation operation $h \mapsto \xi_b^{\St} h$ on $\cH(T)$ at least one of the monomials remains in Equation~\eqref{eqn:sattransform}. Observe that the scalars in front of the monomials in Equation~\eqref{eqn:sattransform} all have the same sign, and that to get the truncated trace on the Steinberg representation we evaluate these monomials at a certain, nonzero point. Thus, the only problem is to see that there is at least one monomial $X$ occurring in $\cS(f_{G, \mu, \alpha})$ and surviving the truncation $X \mapsto \xi_b^{\St} X$. At this point it will be useful to give a graphical interpretation of this truncation process.

A remark on the notation: With $\xi_b^{\St} X$ for $X$ a monomial in the Satake algebra of $T$, we mean the element $\cS_T( \xi_b^{\St} \cS_T^{-1}(X))$ of the Satake algebra of $T$. Below we will use similar conventions for the truncations $\chi_N X$, $\widehat \chi_{N_0 \cap M_b} X$ and $\eta_b X$. 

A \emph{graph} in $\Z^2$ is a sequence of points $\uv_0, \uv_1, \ldots, \uv_r$ with $\uv_{i+1} - \uv_{i} = (1, e)$, where $e$ is an integer. To a monomial $X = X_1^{e_1} X_2^{e_2} \cdots X_n^{e_n} \in \C[X_1^{\pm 1}, X_2^{\pm 1}, \ldots, X_n^{\pm 1}]$, with $e_i \in \Z$ and $\sum_{i=1}^n e_i = s$ we associate the graph $\cG_X$ with points\footnote{To obtain the usual convex picture of the Newton polygon we invert the order of the vector $e_1, \ldots, e_n$.} $\uv_0 := (0, 0)$, $\uv_i := \uv_0 + \lhk i, e_{n} + e_{n-1} + \ldots + e_{n+1-i} \rhk \in \Z^2$, for $i = 1, \ldots, n$. Because the sum $\sum_{i=1}^n e_i$ is equal to $s$, we see that the end point of the graph is $(n, s)$. The function $f_{G\alpha\mu}$ is (up to scalar) the elementary symmetric function of degree $s$ in $n$ variables, thus its monomials correspond precisely to the set of graphs that start at the point $(0, 0)$, have end point $(n, s)$ and satisfy $\uv_{i+1} - \uv_{i} \in \{(1, 0), (1,1)\}$ for all $i$.

To the slopes $\lambda_1 \leq \lambda_2 \leq \cdots \leq \lambda_n$ of the isocrystal $b$ we associate the graph $\cG_b$ with points $\uv_0 := (0, 0)$, $\uv_i := \uv_0 + \lhk i, \lambda_1 + \lambda_2 + \ldots + \lambda_i \rhk \in \Z^2$, for $i = 1, \ldots, n$. 

We may now explain the truncation $X \mapsto \xi_b^{St} X$ in terms of graphs. We have $\xi_b^{\St} X = X$ or $\xi_b^{\St} X = 0$. We claim that we have $\xi_b^{\St} X = X$ if the following conditions hold:
\begin{enumerate}
\item[(\textbf{C1})] We have $\cG_b(n) = \lambda \cG_X(n)$ for some positive scalar $\lambda \in \R_{>0}$;
\item[(\textbf{C2})] For every break point $x \in \Z^2$ of $\cG_b$ the point $x$ also lies on the graph $\lambda \cG_X$;
\item[(\textbf{C3})] Outside the set of breakpoints of $\cG_b$, the graph $\lambda \cG_X$ lies strictly below the graph $\cG_b$.
\end{enumerate}
Thus, in short: $\cG_X$ lies below $\cG_b$ and the set of contact points between the two graphs is precisely the begin point, end point and the set of break points of $\cG_b$. See also our preprint \cite{kret2}, where we use this construction in an analogous context.

\begin{remark} In the claim above we say ``if'' and not ``if and only if''. The conditions (\textbf{C1}), (\textbf{C2}) and (\textbf{C3}) are stronger than the condition $\xi_b^{\St} X = X$. In Lemmas~\ref{st:vLemA}, \ref{st:vLemB} and~\ref{st:vLemC} below we give conditions (\textbf{C1'}), (\textbf{C2'}) and (\textbf{C3'}) which, when taken together, are equivalent to ``$\xi_b^{\St} X = X$''. However (\textbf{C1}, \textbf{C2}, \textbf{C3}) is not equivalent to (\textbf{C1'}, \textbf{C2'}, \textbf{C3'}). If you replace $(\textbf{C3'})$ with the stronger condition
``(\textbf{C3''}): \emph{$\li \cG_{x} = \lambda \cG_b$ for some $\lambda \in \R_{>0}$}''
then you have $(\textbf{C1}, \textbf{C2}, \textbf{C3}) \Longleftrightarrow (\textbf{C1'}, \textbf{C2'}, \textbf{C3''})$. 
\end{remark}

Because the above fact is crucial for the argument, let us prove the claim with complete details. Let $X = (e_1, e_2, \ldots, e_n) \in \Z^n = X_*(A_0)$. We want to express the condition $\xi_b^{\St} X = X$ in terms of $\cG_X$. The Satake transform for the maximal torus $T = (\Res_{F^+/\Qp} \Gm)^n$ is simply
\begin{align}\label{eqn:Explicit_Satake_for_Torus}
\cH_0(T) & \isomto \C[X_1^{\pm 1}, X_2^{\pm 1}, \ldots, X_n^{\pm 1}], \cr
\one_{(p^{-e_1} \cO_{F^+}^\times) \times (p^{-e_2} {\cO_{F^+}^\times}) \times \cdots \times (p^{-e_n} {\cO_{F^+}^\times})} & \longmapsto X^{e_1} X^{e_2} \cdots X^{e_n}.
\end{align}
We have $\xi_b^{\St} = \widehat \chi_{N_0 \cap M_b} \chi_{N_b} \eta_b$. Let $(n_a)$ be the composition of $n$ corresponding to the standard parabolic subgroup $P_b$ of $G$. Let $g = (g_1, \ldots, g_n) \in T$ such that $\chi_{N_b}(g) =1$. Explicitly, this means that
\begin{equation}\label{eqn:satB}
|g_1 g_2 \cdots g_{n_1} |^{1/n_1} < |g_{n_1 + 1} g_{n_1 + 2} \cdots g_{n_1 + n_2}|^{1/n_2} < \ldots < |g_{n_{k-1}} g_{n_{k-1} + 1} \cdots g_n |^{1/n_k}, 
\end{equation}
(cf. \cite[Eq. (1.11)]{kret1}). In terms of the graph $\cG_X$ of $X$ this means the following. We have $X \in \ia_0$ and we have the projection $\pi_{0, P_b}(X)$ of $X$ in $\ia_{P_b}$ (obtained by taking the average under the action of the Weyl group of $M_b$). We write $\li \cG_X$ for the graph of $\pi_{0, P_b}(X) \in \ia_{P_b} \subset \ia_0$. This graph $\li \cG_X$ is obtained from the graph $\cG_X$ as follows. Consider the list of points $p_0 := (0,0)$, $p_1 := (n_1, \cG_X(n_1))$, $p_2 := (n_1 + n_2, \cG_X(n_1 + n_2))$, $\ldots$, $p_k := (n, \cG_X(n))$. Connect, using a straight line, the point $p_0$ with $p_1$, and with another straight line, the point $p_1$ with $p_2$, etc, to obtain the graph $\li \cG_X$ from $\cG_X$. From Equations~\eqref{eqn:Explicit_Satake_for_Torus} and~\eqref{eqn:satB} we get:

\begin{lemma}\label{st:vLemA}
For a monomial $X$ we have $\chi_{N_b} X = X$ if and only if (\textbf{C1'}): The graph $\li \cG_X$ is convex. 
\end{lemma}
(Remark: We have $\chi_{N_b} X = 0$ if condition (\textbf{C1'}) is not satisfied. This remark also applies to Lemmas~\ref{st:vLemB} and~\ref{st:vLemC}.) 

Before discussing the function $\widehat \chi_{N_0 \cap M_b}$, let us first discuss in detail the maximal case, \ie the function $\widehat \chi_{N_0}$ for the group $G$ (cf. \cite[Prop. 1.11]{kret1}). We have $\ia_0 = \R^n$, write $H_1, \ldots, H_n$ for the basis of $\ia_0^*$ dual to the standard basis of $\R^n$. Write $\alpha_i$ for root $H_i - H_{i+1}$ in $\ia_0^*$. We have
\begin{equation}\label{eqn:Fundamental_Roots_GLn}
\varpi_{\alpha_i}^G = \lhk H_1 + H_2 + \ldots + H_i - \frac i n \lhk H_1 + H_2 + \ldots + H_n \rhk\rhk \in \ia_0^{G*}. 
\end{equation}
Thus, for a monomial $X = X_1^{e_1} X_2^{e_2} \cdots X_n^{e_n}$ the condition ``$\langle \varpi^G_{\alpha_i}, X \rangle > 0$'' corresponds to 
\begin{equation}\label{eqn:explicietvergelijking}
e_1 + e_2 + \ldots + e_i > \frac i n \lhk e_1 + e_2 + \ldots + e_n \rhk
\end{equation}
Thus we obtain $\cG_X(n + 1 - i) > \frac i n s$, where $s$ is the degree of $X$, \ie $s = \sum_{i=1}^n e_i$. Demanding that $\langle \varpi_\alpha^G, X \rangle$ is positive for all roots $\alpha$ of $G$, is demanding that the graph $\cG_X$ lies strictly below the straight line connecting the point $(0, 0)$ with the point $(n, \cG_X(n))$. (We get `below' and not `above' due to the inversion ``$e_i \mapsto e_{n+1 -i}$'' that we made in the definition of the graph $\cG_X$. )

We now turn to the function $\widehat \chi_{N_0 \cap M_b}$. The group $M_b$ decomposes into a product of general linear groups, say it corresponds to the composition $(n_a)$ of the integer $n$. Thus, the condition ``$\langle \varpi_\alpha^{M_b}, X \rangle > 0$'' is the condition in Equation~\eqref{eqn:explicietvergelijking} but, then for each of the blocks of $M_b$ individually. The conclusion is:

\begin{lemma}\label{st:vLemB}
For any monomial $X$ we have $\widehat \chi_{N_0 \cap X_b} \cdot X = X$ if and only if (\textbf{C3'}): The graph $\cG_X$ lies below $\li \cG_X$ and the two graphs touch precisely at the points $p_i$. 
\end{lemma}

The condition ``$\eta_b X = X$'' means $\pi_{0, P_b}(\Phi(g))$ equals $\lambda \li \nu_b$ for all $g$ lying in the support of the function $\cS_T^{-1}(X)$ on the group $T$. By the explicit formula for the Satake transform (Equation~\eqref{eqn:Explicit_Satake_for_Torus}), the condition is equivalent to the existence of a permutation $w \in \iS_n$ such that the vector 
$$
\lhk \underset {n_1} {\underbrace {e_{w(1)} + e_{w(2)} + \cdots }}, \underset {n_2} {\underbrace {e_{w(n_1+1)} + e_{w(n_1 + 2)} + \cdots }}, \ldots, \underset {n_k} {\underbrace {e_{w(n_1 + n_2 + \ldots + n_{k-1} + 1)} + \cdots }} \rhk \in \ia_{P_b}, 
$$
is a positive scalar multiple of the vector $\li \nu_b$. Using earlier notations we get:

\begin{lemma}\label{st:vLemC}
For any monomial $X$ we have $\eta_b X = X$ if and only if (\textbf {C3'}): There exists an element $w \in \iS_n$ such that $\li \cG_{w(X)} = \lambda \cG_b$ for some $\lambda \in \R_{>0}$.
\end{lemma}

To prove the claim we show that the group of conditions (\textbf{C1}), (\textbf{C2}) and (\textbf{C3}) implies the group of conditions (\textbf{C1'}), (\textbf{C2'}) and (\textbf{C3'}). 

Thus, assume conditions (\textbf{C1}), (\textbf{C2}) and (\textbf{C3}) are true for the monomial $X$. The parabolic subgroup $P_b$ is contracted by the isocrystal $b$. Thus the set of breakpoints of the polygon $\cG_b$ is equal to the set $q_0 = (0,0)$, $q_1 = (n_1, \cG_b(n_1))$, $q_2 = (n_1 + n_2, \cG_b(n_1 + n_2))$, $\ldots$, $q_k = (n, \cG_b(n))$. By condition (\textbf{C1}) there is a $\lambda \in \R_{>0}$ such that $\cG_b(n) = \lambda \cG_X(n)$.
By conditions (\textbf{C2}) and (\textbf{C3}) the set $\{q_0, \ldots, q_n\}$ is then precisely the set of points where the graph $\lambda \cG_X$ touches the graph $\cG_b$. Taking averages, we get the relation $\li \cG_b = \lambda \li \cG_X$. We have $\cG_b = \li \cG_b$ (because $P_b$ is associated to $b$), and therefore $\cG_b = \lambda \li \cG_X$. Thus condition (\textbf{C3'}) is true for $w = \Id \in \iS_n$. (\textbf{C2'}) is now implied by (\textbf{C2}) and (\textbf{C3}). Finally we prove condition (\textbf{C1'}). We have $\lambda \li \cG_X = \cG_b$, and the graph $\cG_b$ is convex. Thus $\li \cG_X$ is convex. The three conditions (\textbf{C1'}), (\textbf{C2'}) and (\textbf{C3'}) are now verified, and therefore the claim is true. 

The monomials $M$ occurring in $\cS(f_{G, \mu, \alpha})$ corresponds to the set of graphs from $(0,0)$ to $(n,s)$ whose steps consist of diagonal, north-eastward steps, or horizontal, eastward steps. Thus, it suffices to show that there exists a graph satisfying (\textbf{C1}), (\textbf{C2}) and (\textbf{C3}). This is indeed possible (see Figure 2 for the explanation). This completes the proof in case $F^+ \otimes E_{\alpha}$ is a field. 

\begin{figure}
\begin{center}
\input{Figure_pathexist.tex}
\end{center}
\begin{caption}{The dark line is an example of the Newton polygon of an isocrystal $b$ with additional $\Gl_{12}(F^+)$-structure whose slope morphism is $(\tfrac 15, \tfrac 15, \tfrac 15, \tfrac 15, \tfrac 15, \tfrac 12, \tfrac 12, \tfrac 12, \tfrac 12, \tfrac 23, \tfrac 23, \tfrac 23)$. The thin line is a $\xi_b^{\St}$-admissible path. For this Newton polygon there exist precisely two admissible paths. In general one takes the `ordinary' path starting with horizontal steps within the blocks where the Newton polygon is of constant slope, and ending with diagonal steps. 
}
\end{caption}
\end{figure}

We now drop the assumption that the algebra $F^+ \otimes E_{\alpha}$ is a field. By \cite[Prop.~3.3]{kret1} there exists a sufficiently large integer $M \geq 1$ such that for all degrees $\alpha$ divisible by $M$, the function $f_{G\mu \alpha}$ is (up to a scalar) a convolution product of the form $\prod_{i=1}^r f_{n\alpha s_i}$, where $r = [F^+:\qp]$ and $(s_i)$ is a certain given composition of an integer $s$ of length $r$. Any monomial occurring in $\cS(f_{n\alpha s})$ also occurs in the product $\prod_{i=1}^r \cS(f_{n\alpha s_i})$ with a positive coefficient. Thus we may write $\prod_{i=1}^r f_{n\alpha s_i} =f_{n\alpha s} + R \in \cH(G)$ for some function $R \in \cH(G)$, whose Satake transform is a linear combination of monomials, with all coefficients positive. Consequently, to check that the truncated trace of $\prod_{i=1}^r f_{n\alpha s_i}$ on the Steinberg representation is non-zero, it suffices to check that the truncated trace of $f_{n\alpha s}$ on Steinberg is non-zero. This completes the proof.
\end{proof}

\begin{proposition}\label{st:glnonzeroclass}
Let $G$ be an algebraic group over $\qp$ defined as in Equation~\eqref{eqn:Local_group_Unitary}, and assume it is of linear type. Let $b \in B(G, \mu)$ be a $\mu$-admissible isocrystal. Let $m_0$ be the number of indices $i$ such that $\lambda_i = 0$, and let $m_1$ be the number of indices $i$ such that $\lambda_i = 1$. Write $m := n - m_0 - m_1$. Let $\pi_{m_0}$ (resp. $\pi_{m_1}$) be any generic unramified representation of $\Gl_{m_0}(F^+)$ (resp. $\Gl_{m_1}(F^+)$), and $\chi$ an unramified character of $\Gl_m(F^+)$. Let $P$ be the standard parabolic subgroup of $G$ with $3$ blocks, the first of size $m_1$, the second of size $m$ and the last one of size $m_3$. Then for $\alpha$ sufficiently divisible we have 
$$
\Tr \lhk \chi_b^G f_{G,\mu,\alpha}, \Ind_{P}^{G}\lhk \pi_{m_1} \otimes \St_{\Gl_m(F^+)}(\chi) \otimes \pi_{m_0}\rhk \rhk \neq 0.
$$
\end{proposition}
\begin{remark}
We have abused language slightly saying that $P$ has $3$ blocks. We could have $m$, $m_0$ or $m_1$ equal to $0$, in which case $P$ has less than $3$ blocks. If one of the numbers $m, m_0$ or $m_1$ is $0$, then one simply removes the corresponding factor from tensor product $\pi_{m_1} \otimes \St_{\Gl_m(F^+)}(\chi) \otimes \pi_{m_0}$, and one induces from a parabolic subgroup with two blocks (or one block). 
\end{remark}
\begin{proof}[Proof of Proposition~\ref{st:glnonzeroclass}]
By van Dijk's formula for truncated traces \cite[Prop. 1.5]{kret1}, we get a trace on $M$:
\begin{equation}\label{eqn:exttA}
\Tr \lhk \chi_b^G f_{G,\mu,\alpha}^{(P)}, \pi_{m_1} \otimes \St_{\Gl_m(F^+)} \otimes \pi_{m_0} \rhk.
\end{equation}
 By Proposition~\ref{st:constantterms} we have
\begin{align}\label{eqn:exttB}
f_{G, \mu, \alpha}^{(P)} = q^{-\alpha \langle \rho_G - \rho_M, \mu \rangle } \sum_{w \in W_\alpha)/\stab_{W_\alpha)}(\mu) W_{M, \alpha}} f_{M,w(\mu),\alpha} \in \cH_0(M).
\end{align}
The intersection $\Omega_{\li \nu_b}^G \cap M$ is equal to a union $\bigcup \Omega_{w(\li \nu_b)}^M$ with $w$ ranging over the permutations $w \in W$ such that $w(\li \nu_b)$ is $M$-positive. Consequently, if we plug Equation~\eqref{eqn:exttB} into Equation~\eqref{eqn:exttA}, then we get a large sum, call it $(\star)$, of traces of functions $f_{M, w(\mu), \alpha}$ against a representation of the form $\pi_{m_1} \otimes \St_{\Gl_m(F^+)} \otimes \pi_{m_0}$. All the signs are the same in this large sum $(\star)$, therefore it suffices that there is at least one non-zero term. Take $b_M \in B(M)$ the isocrystal whose slope morphism is $\lambda_1 \leq \lambda_2 \leq \cdots \leq \lambda_n$ in the $M$-positive chamber of $\ia_0$. Then $b_M$ has only slopes $0$ on the first block of $M$ and only slopes $1$ on the third block, and all its slopes $\neq 0, 1$ are in the second block. The trace $\Tr(\chi_{b_M}^M f_{M,\mu,\alpha}, \pi_{m_1} \otimes \St_{\Gl_m(F^+)} \otimes \pi_{m_0})$ occurs as a term in the expression $( \star )$. By Lemma~\ref{st:prenonzerolemma} and Proposition~\ref{st:ACclassicalgroups} this term is non-zero. This completes the proof.
\end{proof}

We now establish the cases where the group is an unramified unitary group over $F^+$ (unitary type, cf. Equation~\eqref{eqn:Local_group_Unitary}). 

\begin{lemma}\label{st:lemmasteinunitary}
Let $G$ be an algebraic group over $\qp$ defined as in Equation~\eqref{eqn:Local_group_Unitary}, and assume it is of unitary type. Let $b \in B(G, \mu)$ be an $\mu$-admissible isocrystal whose slope morphism $\li \nu_b \in \ia_0$ has no coordinate equal to $0$ and no coordinate equal to $1$. Then, for $\alpha$ sufficiently divisible, the trace $\Tr(\chi_b^{G} f_{G, \mu, \alpha}, \St_{G(\qp)})$ is non-zero. 
\end{lemma}
\begin{proof} 
We use the explicit description $f_{G, \mu, \alpha} = \Psi f_{G^+, \mu, (\alpha/2)}$ of the Kottwitz function from Equation~\eqref{eqn:Base_Change_Map_Unitary_Group}. Assume the algebra $F^+ \otimes E_\alpha$ is a field; then the base change mapping from $G(F_{\alpha}^+) \to G(F_2)$ is given by $X_i \mapsto X_i^{\alpha/2}$ on the Satake algebras. Over $F_\alpha^+$, the Weyl group $W_\alpha$ is equal to $\iS_n$ with its natural action on $\R^n$. The formula for the base change mapping $\Psi$ from Equation~\eqref{eqn:Base_Change_Map_Unitary_Group} also makes sense over the Satake algebras of the maximal split tori, \ie we have a map $\Psi$ from the algebra $\C[X_1^{\pm 1}, \ldots, X_n^{\pm 1}]$ to $\C[X_1^{\pm 1}, \ldots, X_n^{\pm 1}]$. The monomials occurring in $f_{G, \mu, \alpha}$ are those monomials of the form $\Psi[w(\mu)]$ where $w$ is some element of $\iS_n$. The Weyl group translates $[w(\mu)]$ of $[\mu]$ correspond to all paths from $(0,0)$ to $(n, s)$, and the monomials of the form $\Psi [w(\mu)] = [w(\mu)] + [\theta(w\mu)]$ correspond to all paths from $(0,0)$ to $(n, 0)$ staying below the horizontal line with equation $y = s$, and above the horizontal line with equation $y = -s$. The truncation $\chi_b^{G(\qp)} \Psi[w(\mu)]$ is non-zero if the path $\cG$ of $\Psi[w(\mu)]$ lies below $\cG_b$ and the set of contact points between the two graphs is precisely the initial point, end point and the set of break points of $\cG_b$. This is the same condition as had for the general linear group (see above Equation~\eqref{eqn:Explicit_Satake_for_Torus}) because the root systems are the same. Such graphs exist in case $b$ has no slopes equal to $-1, 0$ or $1$ (draw a picture). Consequently $\chi_b^{G(\qp)} f_{G, \mu, \alpha} \neq 0$, and then also $\Tr(\chi_b^{G(\qp)} f_{G, \mu, \alpha}, \St_G) \neq 0$ by Proposition~\ref{st:truncatedtraceonsteinberg}. 

Forget the assumption that $F^+ \otimes E_{\alpha}$ is a field. We proceed just as we did for the general linear group (cf. Lemma~\ref{st:ACclassicalgroups}), we write $f_{G, \mu, \alpha} = A + R$, where $R$ is a function whose Satake transform is a linear combination of monomials in the Satake algebra with all coefficients positive, and $A$ is a function for which we already know that its truncated trace on the Steinberg representation does not vanish. This completes the proof.
\end{proof}

\begin{proposition}\label{st:unitnonzeroclass}
Let $G$ be an algebraic group over $\qp$ defined as in Equation~\eqref{eqn:Local_group_Unitary}, and assume it is of unitary type. Let $b \in B(G, \mu)$ be an isocrystal with slopes $\lambda_1 \leq \lambda_2 \leq \cdots \leq \lambda_n$ (cf. the discussion below Proposition~\ref{st:slopeproposition}). Let $n = m_1 + m_2 + m_3$ be the composition of $n$ such that the first block of $m_1$ slopes $\lambda_i$ satisfy $\lambda_i = -1$, the second block of slopes $\lambda_i$ satisfy $-1 < \lambda_i < 1$ and is of size $m_2$, the third block of slopes $\lambda_i$ satisfy $\lambda_i = 1$ and is of size $m_3$. We have $m_1 = m_3$. Let $P = MN$ be the standard parabolic subgroup of $G$ corresponding to this composition of $n$, thus $M$ is a product of two groups, $M = M_1 \times M_2$, where $M_1 = \Gl_{m_1}(F^+)$ is a general linear group and $M_2$ is an unramified unitary group. For $\alpha$ sufficiently divisible the trace $\Tr(\chi_b^{G(F^+)} f_{G, \mu, \alpha}, \bullet)$ against the representation $\Ind_{P(F^+)}^{G(F^+)} (\pi_{m_1} \otimes \St_{m_2}(\chi))$ is non-zero if $\pi_{m_1}$ is an unramified generic representation and $\chi$ an unramified character of $\Gl_{m_2}(F^+)$. 
\end{proposition}
\begin{remark}
The group $M_1$ could be trivial. This happens in case $-1 < \lambda < 1$ for all indices $i$. When $M_1$ is trivial, the considered representation is simply an unramified twist of the Steinberg representation. 
\end{remark}
\begin{proof}
The proof is the same as the proof in case of the general linear group (cf. Proposition~\ref{st:glnonzeroclass}): one easily reduces the statement to Lemma~\ref{st:lemmasteinunitary}. 
\end{proof}

Let now $G/\Q$ be an unitary group of similitudes arising from a Shimura datum of PEL-type, and let $G_1 \subset G$ be the kernel of the factor of similitudes. The group $G_1$ is defined over a totally real field $F^+$, and defined with respect to a quadratic extension $F$ of $F^+$, which is a CM field. Let $A_0 \subset G$ be a maximally split torus, then we may write $A_0 = \Gm \times A'_0$ (\emph{not} a direct product), where $A'_0 \subset G_1$ be the maximally split torus of $G_1$ defined by $G_1 \cap A_0$. At $p$ we have a decomposition of $F^+ \otimes \qp$ into a product of fields $F^+_\wp$, where $\wp$ ranges over the primes above $p$. Let $p$ be a prime number where $G$ is unramified. The group $G_{1, \qp}$ is of the form $G_{1, \qp} \cong \prod_{\wp} \Res_{F_{\wp}^+/\qp} G_{1,\wp}$, where the group $G_{1,\wp}$ is either an unramified unitary group over $F_\wp^+$, or the general linear group. In the first case we call the $F^+$-prime $\wp$ \emph{unitary} and in the second case we call the prime \emph{linear}. 

Consider an isocrystal $b \in B(G)$. To $b$ we may associate its slope morphism $\li \nu_b \in \ia_0$. Let $A_{0,\wp}' \subset G_{1,\wp}$ be the $\wp$-th component of $A_0'$; it is a split maximal torus in $G_{1,\wp}$, and write $\ia_0(\wp) := X_*(A_{0,\wp}')$. The space $\ia_0$ decomposes along the split center and the $F^+$-primes $\wp$ above $p$: $\ia_0 = \R \times \prod_{\wp} \ia_0(\wp)$. We speak for each $\wp$ of the $\wp$-component $\li \nu_{b, \wp}$ of $\li \nu_b$. In case $\wp$ is linear, the Proposition~\ref{st:glnonzeroclass} gives us a class of representations $\pi_\wp'$ of $G_{1,\wp}(\qp)$ such that the $b_\wp$-truncated trace on $\pi_\wp'$ does not vanish. In case $\wp$ is unitary, we get such a class $\pi_\wp'$ from Proposition~\ref{st:unitnonzeroclass}. Let $\pi'$ be the representation of $G_1(\qp)$ obtained from the factors $\pi_\wp$ by taking the tensor product.

\begin{definition}
We write $\iR_1(b)$ for the just constructed class of $G_1(\qp)$-representations $\pi'$. 
\end{definition}

\begin{remark}
The set of representations $\iR_1(b)$ has positive Plancherel measure in the set of $G_1(\qp)$ representations, and the $b$-truncated trace of the Kottwitz function on these representations does not vanish by construction.
\end{remark}

We now extend the class $\iR_1(b)$ to a class of $G(\qp)$-representations, as follows:

\begin{definition}\label{classSb}
Let $\pi \in \iR_1(b)$. Then $\pi$ is an $H(\qp)$-representation; let $\omega_\pi$ be its central character, thus $\omega_\pi$ is a character of $Z_1(\qp)$. Assume $\chi$ is a character of $Z(\qp)$ extending $\omega_\pi$. Then we may extend the representation $\pi$ to a representation $\pi\chi$ of the group $H(\qp)Z(\qp)$. We define $\iR_1(b)'$ to be the set of $H(\qp)Z(\qp)$-representations of the form $\pi\chi$. Not all the inductions $\Ind_{H(\qp)}^{G(\qp)} (\pi \chi)$ have to be irreducible, we ignore the reducible ones. We define $\iR(b)$ to be the set of representations $\Pi$ isomorphic to an irreducible induction $\Ind_{H(\qp)Z(\qp)}^{G(\qp)}(\pi\chi)$ with $\pi\chi \in \iR_1(b)'$.
\end{definition}

The required non-vanishing property of the representations in $\iR_b$ will be shown in the next section. 

\subsection{Local extension}
We need to extend from $G_1(\qp)$ to the group $G(\qp)$. Let $Z$ be the center of the group $G$. Consider the morphism of algebraic groups $\psi \colon G_{1, \qp} \times Z_{\qp} \surjects G_{\qp}$; the group $\Ker(\psi)$ is the center $Z_1$ of the group $G_1$, so
\begin{equation}\label{eqn:kerpsi}
\Ker(\psi) = \prod_\wp \begin{cases} \Gm & \wp \textup{ is linear} \cr
U_1^* & \wp \textup{ is unitary,} 
\end{cases}
\end{equation}
where $U_1^*$ is the unramified non-split form of $\Gm$ over $F_{\wp}^+$. Over $\Q$, $Z$ is defined by $Z(\Q) = \{x \in F^\times | N_{F/F^+}(x) \in \Q^\times \}$. Using Equation~\eqref{eqn:kerpsi}, the long exact sequence for Galois cohomology and Shapiro's lemma, the group $G(\qp)/G_1(\qp)Z(\qp)$ maps injectively into the group $(\Z/2\Z)^t$, where $t$ is the number of unitary places of $F^+$ above $p$. 

Write $\mu' \in X_*(T)$ for the cocharacter of the maximal torus $(T \cap G_1) \cap Z$ of $G_1 \times Z$ obtained from $\mu$ via restriction. Let $f_{G_1 \times Z, \mu', \alpha}$ be the corresponding function of Kottwitz on the group $G_1(\qp) \times Z(\qp)$. 
Furthermore we write $\chi_b^{G_1 \times Z}$ for the characteristic function on $G_1(\qp) \times Z(\qp)$ of elements $(g, z)$ such that we have $\chi_b^{G(\qp)}(gz) = 1$. 
We prove the following statement:

\begin{proposition}\label{st:atp}
Fix a representation $\pi_0$ of $G_1(\qp$). Let $\Pi$ be a smooth irreducible representation of $G(\qp)$ containing the representation $\pi_0$ of $G_1(\qp)$ upon restriction to $G_1(\qp) \times Z(\qp)$. Assume the central character of $\Pi$ is of finite order. Then, for all sufficiently divisible $\alpha$, we have 
$$
\Tr(\chi_b^{G(\qp)} f_{G, \mu, \alpha}, \Pi) = t(\Pi) \Tr(\chi_b^{G_1 \times Z} f_{G_1 \times Z, \mu, \alpha}, \pi_0)
$$
where $t(\Pi)$ is a positive real number. 
\end{proposition}

Before proving Proposition~\ref{st:atp} we first establish some technical results. We fix smooth models of $G, G_1, Z$, etc. over $\Zp$ (and use the same letter for them). We have the exact sequence $Z_1 \injects Z \times G_1 \surjects G$, so the cokernel of $Z(\qp)G_1(\qp)$ in $G(\qp)$ is a subgroup of $H^1(\qp, Z_1) \cong (\Z/2\Z)^t$, where $t$ is the number of unitary places. 

\begin{lemma}\label{st:thishastobechanged}
The mapping $G_1(\zp) \times Z(\zp) \to G(\zp)$ is surjective.
\end{lemma}
\begin{proof}
We have an exact sequence $Z_1 \injects G_1 \times Z \surjects G$ of algebraic groups over $\spec(\zp)$. Thus we get $Z_1(\fp) \injects G_1(\fp) \times Z(\fp) \to G(\fp) \to H^1(\F_p, Z_1)$. The group $Z_1$ is a torus and therefore connected. By Lang's theorem we obtain $H^1(\F_p, Z_1) = 1$. Thus the mapping $G_1(\fp) \times Z(\fp) \to G(\fp)$ is surjective. By Hensel's lemma the mapping $G_1(\zp) \times Z(\zp) \to G(\zp)$ is then also surjective. 
\end{proof}

\begin{lemma}\label{st:supportlemma}
The function of Kottwitz $f_{G, \mu, \alpha}$ has support on the subset $Z(\qp)G_1(\qp) \subset G(\qp)$. 
\end{lemma}
\begin{proof}
Define $\chi$ on $G(\qp)$ to be the characteristic function of the subset $Z(\qp) G_1(\qp) \subset G(\qp)$. The mapping $Z \times G_1 \to G$ is surjective on $\Zp$-points, and therefore $\chi$ is spherical. The functions $\chi f_{G, \mu, \alpha}^G$ and $f_{G, \mu, \alpha}$ are then both spherical functions and to show that they are equal it suffices to show that their Satake transforms agree (the Satake transform is injective). We have $\cS(\chi f_{G, \mu, \alpha}) = \chi|_{A(\qp)} \cS(f_{G, \mu, \alpha})$, where $A$ is a maximal split torus of $G$, $\chi|_{A(\qp)}$ is the characteristic function of the subset $Z(\qp)G_1(\qp) \cap A(\qp) \subset A(\qp)$. In fact, $Z(\qp)G_1(\qp) \cap A(\qp) = A(\qp)$. Thus $\chi|_{A(\qp)} \cS(f_{G, \mu, \alpha}) = \cS(f_{G, \mu, \alpha})$, and the functions $\chi f_{G, \mu, \alpha}$ and $f_{G, \mu, \alpha}$ have the same Satake transform. This completes the proof of the lemma.
\end{proof}

We now turn to the proof of Proposition~\ref{st:atp}.

\begin{proof}[Proof of Proposition~\ref{st:atp}] 
By Clifford theory \cite[thm 2.40]{chunhui} the representation $\Pi$ restricted to $G_1(\qp)Z(\qp)$ is a finite direct sum of irreducible representations $\pi_i$, where $\pi_i$ satisfies $\pi_i(g) = \pi_0(x_i g x_i^{-1})$ for some $x_i$ not depending on $g$. (In this finite direct sum multiplicities may occur.) As characters on $G_1(\qp)Z(\qp)$ we may write $\theta_\Pi = \sum_{i=1}^t \theta_{\pi_i} \omega_i$, where $\theta_{\pi_i}$ is the Harish-Chandra character of $\pi_i$, viewed as a $G_1(\qp)$-representation, and $\omega_i$ is the central character of $\pi_i$. Using Lemma~\ref{st:supportlemma} we may now compute:
\begin{align}
\Tr(\chi_b^{G(\qp)} f_{G, \mu, \alpha}, \Pi) & = \int_{Z(\Qp)G_1(\qp)} \chi_b^{G(\qp)} f_{G, \mu, \alpha} \theta_\Pi \dd g \cr
& = \sum_{i = 1}^t \int_{Z(\qp)G_1(\qp)} \chi_b^{G(\qp)} f_{G, \mu, \alpha} \theta_{\pi_i} \omega_i \dd g \cr
& = \sum_{i=1}^t \int_{Z(\qp)G_1(\qp)} \chi_b^{G(\qp)} f_{G, \mu, \alpha}^{x_i^{-1}} \theta_{\pi_0} \omega_0 \dd g, 
\end{align}
where $f_{G, \mu, \alpha}^{x_i^{-1}}$ is the conjugate of $f_{G, \mu, \alpha}$ by $x_i^{-1}$. However, the function of Kottwitz \emph{is} stable under the action of the Weyl group of $G$. Therefore $f_{G, \mu, \alpha}^{x_i^{-1}} = f_{G, \mu, \alpha}$. We get the expression:
$$
t \int_{Z(\qp)G_1(\qp)} \chi_b^{G(\qp)} f_{G, \mu, \alpha} \theta_{\pi_0} \omega_0 \dd g.
$$
On the other hand we have 
$$
0 \neq \Tr(\chi_b^{Z \times G_1} f_{Z \times G_1, \mu', \alpha}, \pi_0) = \int_{Z(\qp) \times G_1(\qp)} \chi_b^{Z \times G_1} f_{Z \times G_1, \mu', \alpha} [\theta_{\pi_0} \times \omega_0] \dd g. 
$$
We compute the right hand side:
\begin{align}\label{eqn:rightside}
\int_{\frac { Z(\qp) \times G_1(\qp) } {Z_1(\qp)}} & \int_{Z_1(\qp)} (\chi_b^{Z \times G_1} f_{Z \times G_1, \mu', \alpha} [\theta_{\pi_0} \times \omega_0)](z z_1, h z_1) \dd z_1 \frac {\dd (z, h)} {\dd z_1} \cr
 & =\int_{\frac { Z(\qp) \times G_1(\qp) } {Z_1(\qp)}} \chi_b^{Z \times G_1} \int_{Z_1(\qp)} f_{Z \times G_1, \mu', \alpha}(zz_1, hz_1) \dd z_1 (\theta_{\pi_0} \omega_0)(z, h) \frac {\dd (z, h)}{\dd z_1}. 
\end{align}
We claim that
\begin{equation}\label{eqn:integralclaim}
\int_{Z_1(\qp)} f_{Z \times G_1, \mu', \alpha}(zz_1, hz_1) \dd z_1 = f_{G, \mu, \alpha}(z, h).
\end{equation}
The map $Z \times G_1 \to G$ is surjective on $\zp$-points, and therefore the function 
$$
\int_{Z_1} f_{Z \times G_1, \mu', \alpha}(zz_1, hz_1) \dd z_1 
$$
is $G(\zp)$-spherical. Therefore, to show that Equation~\eqref{eqn:integralclaim} is true, it suffices to show that the Satake transforms of these functions agree. 

We compute the Satake transform of the left hand side:
\begin{align*}
\delta_{P_0}^{-1} \int_{N_0(\qp)} \int_{Z_1(\qp)} & f_{Z \times G_1, \mu', \alpha} (z z_1 n_0, h z_1 n_0) \dd z_1 \dd n_0 \cr 
&= \delta_{P_0}^{-1} \int_{Z_1(\qp)} \int_{N_0(\qp)} f_{Z \times G_1, \mu', \alpha}(z z_1 n_0, h z_1 n_0) \dd n_0 \dd z_1 \cr
& = \int_{Z_1(\qp)} \delta_{P_0}^{-1} \int_{N_0(\qp)} f_{Z \times G_1, \mu', \alpha}(z z_1 n_0, h z_1 n_0) \dd n_0 \dd z_1 \cr
& = \int_{Z_1(\qp)} \lhk f_{Z \times G_1, \mu', \alpha} \rhk^{(P_0)}(z z_1, h z_1) \dd z_1
\end{align*}
By Definition~\ref{def:kottwitzfunction} the last expression is equal to $f_{G, \mu, \alpha}^{(P_0)}(z, h)$. This proves Equation~\eqref{eqn:integralclaim}. We may continue with Equation~\eqref{eqn:rightside} to obtain
$$
\int_{\frac {Z(\qp) \times G_1(\qp)} {Z_1(\qp)}} \chi_b^{Z \times G_1} f_{G, \mu, \alpha} \theta_{\pi_0} \omega_0 \frac {\dd (z, h)}{\dd z_1}. 
$$
Now $\omega_0$ is of finite order by assumption, and the function $f_{G\mu\alpha}$ restricted to $\qp^\times \cong A(\qp) \subset Z(\qp)$ is the characteristic function of $p^{-\alpha}\Zp^\times$. For $\alpha$ sufficiently divisible this is then, up to normalization of Haar measures, just the trace $\Tr(\chi_b^{Z \times G_1} f_{Z \times G_1, \mu', \alpha}, \pi_0)$. This proves that $\Tr(\chi_b^G f_{G, \mu, \alpha}, \Pi)$ and $\Tr(\chi_b^{Z \times G_1} f_{Z \times G_1, \mu', \alpha}, \pi_0)$ differ by a positive, non-zero, scalar. The proof of the theorem is now complete.
\end{proof}

\subsection{The isolation argument}
Let $\Sh_K$ be a Shimura variety of PEL-type of type (A), and let $G$ be the corresponding unitary group of similitudes over $\Q$. We write $E$ for the reflex field and we let $p$ be a prime of good reduction\footnote{Here `\emph{good reduction}' is in the sense of Kottwitz \cite[\S 6]{MR1124982}; in particular $K$ decomposes into a product $K = K_p K^p$ with $K_p \subset G(\qp)$ hyperspecial.}. Let $b \in B(G_{\qp}, \mu)$ be an admissible isocrystal. Let $\p$ be a prime of the reflex field $E$ above $p$. Let $\fq$ be the residue field of $E$ at $\p$. Let $\Sh_{K,\p}^b$ be the corresponding Newton stratum of $\Sh_{K, \p}$, a locally closed subvariety of $\Sh_{K, \p}$ over $\fq$ \cite{MR1411570}.

Let $\alpha$ be a positive integer. We fix an embedding $E_\p \subset \lqp$ and we write $E_{\p, \alpha}$ for the extension of the field $E_{\p}$ of degree $\alpha$ inside $\lqp$. 

\begin{theorem}[Viehmann-Wedhorn]\label{st:mainthmA}
The variety $\Sh_{K, \p}^b$ is not empty. 
\end{theorem}
\begin{remark}
In the statement of the above theorem we have not been precise about the form of the compact open subgroup $K \subset G(\Af)$. Note however that for any pair $(K, K')$ of compact open subgroups, hyperspecial at $p$, we have the finite \'etale morphisms $\Sh_{K} \leftarrow \Sh_{K \cap K'} \rightarrow \Sh_{K'}$ respecting the Newton stratification modulo $\p$. Therefore, showing the Newton stratum is non-empty for one $K$ is equivalent to showing it is non-empty for all $K$. 
\end{remark}
\begin{proof}
Fix a sufficiently divisible and even integer $\alpha$ such that the conclusion of Proposition~\ref{st:atp} is true. We start with the formula of Kottwitz. We write $\phi_\alpha$ for the function $\phi_{G, \mu, \alpha}$ from the previous section\footnote{Where the notation $E_\alpha$ from that section should be replaced with $E_{\p, \alpha}$, and similarly $F^+$ of that section should be replaced by the algebra $F^+ \otimes \Qp = \prod_{\wp} F^+_\wp$, where $\wp$ ranges over the places above $p$.} on $G(E_{\p, \alpha})$. Similarly $f_\alpha := f_{G, \mu, \alpha}$. We pick a prime $\ell \neq p$ and fix an isomorphism $\lql \cong \C$ (and suppress it from all notations). Let $\xi$ be an irreducible complex (algebraic) representation of $G$, and write $\cL$ for the corresponding $\ell$-adic local system on the Shimura tower. Then the Kottwitz formula states:
\begin{equation}\label{eqn:Kottwitz_Formula_Original}
\sum_{x' \in \Fix_{f^p \times \Phi_\p^\alpha}^b(\lfq)} \Tr(f^p \times \Phi_\p^{\alpha}, \iota^* (\cL)_x) = |\Ker^1(\Q, G)| \sum_{(\gamma_0; \gamma, \delta)} c(\gamma_0; \gamma, \delta) O_\gamma(f^{\infty p}) \TO_\delta(\phi_\alpha) \Tr \xi_\C(\gamma_0),
\end{equation}
where $\Fix_{f^p \times \Phi_\p^\alpha(\lfq)}^b$ is the set of fixed points of the Hecke correspondence $f^p \times \Phi_\p^\alpha$ acting on $\Sh_{K, \fq}^b$, and where the sum ranges of the Kottwitz triples $(\gamma_0; \gamma, \delta)$ with the additional condition that the isocrystal defined by $\delta$ is equal to $b$. In Equation~\eqref{eqn:Kottwitz_Formula_Original} the map $\iota$ is the embedding of $\Sh_{K, \Fq}^b$ into $\Sh_{K, \fq}$. 

We may rewrite the right hand side of Equation~\eqref{eqn:Kottwitz_Formula_Original} as
\begin{equation}\label{eqn:Kottwitz_Formula_Rewrite}
|\Ker^1(\Q, G)| \sum_{(\gamma_0; \gamma, \delta)} c(\gamma_0; \gamma, \delta) \cdot O_\gamma(f^{\infty p}) \TO_\delta(\chi_{\sigma b}^{G(E_{\p, \alpha})} \phi_\alpha) \Tr \xi_\C(\gamma_0),
\end{equation}
where now the sum ranges over \emph{all} Kottwitz triples and where $\chi_{\sigma b}^{G(E_{\p, \alpha})}$ is the characteristic function on $G(E_{\p, \alpha})$ such for each element $\delta \in G(E_{\p, \alpha})$ we have $\chi_{\sigma b}^{G(E_{\p, \alpha})}(\delta) = 1$ if and only if the conjugacy class $\gamma = \cN(\delta)$ satisfies $\Phi(\gamma) = \lambda \li \nu$ for some positive real number $\lambda \in \R^\times_{>0}$. Assume the triple $(\gamma_0; \gamma, \delta)$ is such that the corresponding term $c(\gamma_0; \gamma, \delta) O_\gamma(f^{\infty p}) \TO_\delta(\chi_{\sigma b}^{G(E_{\p, \alpha})} \phi_\alpha) \Tr \xi_\C(\gamma_0)$ is non-zero. Then, by the proof of Kottwitz \cite{MR1044820}, the triple $(\gamma_0; \gamma, \delta)$ arises from some virtual Abelian variety with additional PEL-type structures. In particular the isocrystal defined by $\delta$ lies in the subset $B(G_{\qp}, \mu) \subset B(G_{\qp})$. Thus its end point is determined. We have $\gamma = \cN(\delta)$ and $\Phi(\gamma) = \lambda \li \nu_b$ for some $\lambda$ (Proposition~\ref{st:slopeproposition}). Therefore the isocrystal defined by $\delta$ must be equal to $b$. Thus the above sum precisely counts Abelian varieties with additional PEL type structures over $\F_{q^\alpha}$ such that their isocrystal equals $b$. 

We show that the sum in Equation~\eqref{eqn:Kottwitz_Formula_Rewrite} is non-zero. Let $\cE$ be the (finite) set of endoscopic groups $H$ associated to $G$ and unramified at all places where the data $(G, K)$ are unramified. By the stabilization argument of Kottwitz \cite{MR1163241}, the expression in Equation~\eqref{eqn:Kottwitz_Formula_Rewrite} is equal to the stable sum 
\begin{equation}\label{eqn:Kottwitz_Formula_Stable_Geometric}
 \sum_{\cE} \iota(G, H) \cdot \ST_e^*((\chi_b^G f_\alpha)^H), 
\end{equation}
where $(\chi_b^G f_\alpha)^H$ are the transferred functions, whose existence is guaranteed by the fundamental lemma, the $*$ in $\ST_e^*$ means that one only considers stable conjugacy classes satisfying a certain regularity condition (which is empty in case $H$ is a maximal endoscopic group), and finally $\iota(G, H)$ is a constant depending on the endoscopic group (cf. [\textit{loc. cit.}] for the definition). 

We consider only functions such that the transfer $(\chi_b^G f_\alpha)^H$ vanishes for proper endoscopic groups, and therefore we may ignore the regularity condition\footnote{In fact, due to the form of the function $f_\infty$ we have $\ST_e^* = \ST_e$, see \cite[thm 6.2.1]{MR2567740} or \cite[(2.5)]{clozelpurity}.}. Thus, Equation~\eqref{eqn:Kottwitz_Formula_Stable_Geometric} simplifies for such functions and gives the equation:
\begin{equation}\label{eqn:Kottwitz_Formula_Truncated}
\sum_{x' \in \Fix_{f^p \times \Phi_\p^\alpha}^b(\lfq)} \Tr(f^p \times \Phi_\p^{\alpha}, \iota^* (\cL)_x) = \sum_{\cE} \iota(G, H) \ST_e((\chi_b^G f_\alpha)^H).
\end{equation}

Visibly, if the left hand side of Equation~\eqref{eqn:Kottwitz_Formula_Truncated} is non-zero for some Hecke operator $f^p$, then the variety $\Sh_{K, \fq}^b$ is non-empty. We will show that the right hand side of Kottwitz's formula does not vanish for some choice of $K^p$ and some choice of $f^p$. 

We write $G_0^*, G_1^*, G^*$ for the quasi-split inner forms of $G_0$, $G_1$, and $G$ respectively (we remind the reader that $G_0$ is defined over $F^+$ and that $G_1 = \Res_{F^+/\Q} G_0$). The group $G^*$ is the maximal endoscopic group of $G$. Let $\{x_1, x_2, \ldots, x_d\}$ be the set of prime numbers such that the group $G_{\Q_{x_i}}$ is ramified. For $v$ a prime number with $v \notin \{x_1, x_2, \ldots, x_d\}$ the local group $G_{\Q_v}$ is quasi-split, and therefore we may (and do) identify it with the group $G_{\Q_v}^*$. Below we will transfer functions from the group $G(\A)$ to the group $G^*(\A)$; at the places $v$ with $v \notin \{x_1, x_2, \ldots, x_d, \infty\}$ we have $G(\Q_v) = G^*(\Q_v)$ and using this identification we may (and do) take $(h_v)^{G^*(\Q_v)} = h_v$ for any $h_v \in \cH(G(\Q_v))$. 

\begin{remark}
To help the reader understand what we do below at the places $x_i$ (and \emph{why} we do this), let us interrupt this proof with a general remark on the fundamental lemma. It is important to realise that if $v = x_i$ is one of the bad places, then the fundamental lemma guarantees the \emph{existence} of the transferred function $h_v \rightsquigarrow (h_v)^{G^*(\Q_v)}$; however, in its current state, the fundamental lemma does \underline{not} give an explicit description of a transferred function $(h_v)^{G^*(\Q_v)}$. The fundamental lemma only gives explicit transfer in case the group is unramified and the level is hyperspecial. In our case the transferred function $(h_v)^{G^*(\Q_v)}$ is not explicit, and this could introduce signs and cancellations we cannot control. This makes it hard to show that expressions such as the one in Equation~\eqref{eqn:Arthur_Formula} do not vanish. In the argument below we solve the issue by taking $h_v$ to be a pseudocoefficient of the Steinberg representation. For these functions an explicit transfer is known (the transfer is again a pseudocoefficient of the Steinberg representation) and therefore we will be able to control the signs and avoid cancellations. 
\end{remark}

We are going to construct an automorphic representation $\Pi_0$ of $G^*$ with particularly nice properties. From this point onward we take $\xi$ to be a fixed, sufficiently regular complex representation (in the sense of \cite[Hyp.~(1.2.3)]{MR2856383}). We also assume that $\xi$ defines a coefficient system of weight $0$ (cf. \cite{clozelpurity}), and, even better, $\xi$ is trivial on the center of $G^*$. Fix three additional, different, prime numbers $p_1, p_2, p_3$ ($\neq p$) such that the group $G_{\Q_{p_i}}$ is split for $i=1,2,3$. Let $\Pi_{0, p_1}$ be a cuspidal representation of the group $G(\Q_{p_1}) = G^*(\Q_{p_1})$. Let $A(\R)^+$ be the topological neutral component of the set of real points of the split center $A$ of $G$. We apply a theorem of Clozel and Shin \cite{shinplancherel, MR818353} to find an automorphic representation $\Pi_0 \subset L^2_0(G^*(\Q)A(\R)^+ \backslash G^*(\A))$ of $G^*(\A)$ with:
\begin{enumerate}
\item $\Pi_{0,\infty}$ is in the discrete series and is $\xi$-cohomological;
\item $\Pi_{0, p}$ lies in the class $\iR(b)$ (cf. Definition \ref{classSb});
\item $\Pi_{0, p_1}$ lies in the \emph{inertial orbit}\footnote{For the definition of inertial orbit, see \cite[V.2.7]{MR2567785}.} $\cI(\Pi_{0, p_1})$ of $\Pi_{0, p_1}$ at $p_1$;
\item $\Pi_{0, p_2}$ is isomorphic to the Steinberg representation (up to an unramified twist of finite order);
\item $\Pi_{0, x_i}$ is isomorphic to an unramified twist (of finite order) of the Steinberg representation of $G(\Q_{x_i})$ (for $i=1, 2, \ldots, d$);
\item $\Pi_{0,v}$ is unramified for all primes $v \notin \{p, p_1, p_2, p_3, x_1, x_2, \ldots, x_d\}$;
\item The central character of $\Pi_0$ has finite order.
\end{enumerate}
Because the component at $p_1$ of $\Pi_{0}$ is cuspidal, the representation $\Pi_0$ is a cuspidal automorphic representation. The point (7) is possible because of the condition on the weight of $\xi$. 

We now choose the group $K \subset G(\Af)$, and we will also choose a compact open group $K^*$ in $G^*(\Af)$. Write $S = \{p, p_1, p_2, p_3\}$. Write $S' = \{p, p_1, p_2, p_3, x_1, x_2, \ldots, x_d\}$ for the union of $S$ with the set of all places where the group $G$ is ramified. 

The compact open group $K \subset G(\Af)$ is a (any) group with the following properties: 
\begin{enumerate}
\item $K$ is a product $\prod_v K_v \subset G(\Af)$ of compact open groups;
\item for all $v \notin S'$ the group $K_v$ is hyperspecial;
\item $K_p$ is hyperspecial;
\item $K_{p_3}$ is sufficiently small so that $\Sh_{K}$ is smooth and $(\Pi_{0,p_3})^{K_{p_3}} \neq 0$;
\item $K_{x_i}$ is sufficiently small so that the function $f_{x_i}$ is $K_{x_i}$-spherical;
\item for all $v \notin \{x_1, x_2, \ldots, x_d\}$ the space $(\Pi_{0, v})^{K_v}$ is non-zero. 
\end{enumerate}
The group $K^* \subset G^*(\Af)$ is a (any) group with the following properties:
\begin{enumerate}
\item $K^*$ is a product $\prod_v K_v^* \subset G^*(\Af)$ of compact open groups;
\item for any prime $v \notin \{x_1, \ldots, x_d\}$ we have $K_v^* = K_v \subset G(\Q_v) = G^*(\Q_v)$;
\item for all $i \in \{1, 2, \ldots, d\}$ we have $(\Pi_{0, x_i})^{K_{x_i}} \neq 0$;
\end{enumerate}
We now choose the Hecke function $f \in \cH(G(\Af))$. Consider the function $f^{p\infty} \in \cH(G(\Af))$ of the form 
\begin{equation}\label{eqn:The_Hecke_Operator}
f^{p\infty} := f_{p_1} \otimes f_{p_2} \otimes f_{p_3} \otimes f_{x_1} \otimes f_{x_2} \otimes \cdots \otimes f_{x_d} \otimes f^{S'},
\end{equation}
 where 
\begin{itemize}
\item $f_{p_1}$ is a pseudo-coefficient on $G(\Q_{p_1})$ of the representation $\Pi_{p_1}$;
\item $f_{p_2}$ is a sign\footnote{See the remark below Proposition~\ref{prop:stein_coefficients} for the value of this sign.} times a pseudo-coefficient of the Steinberg representation of $G(\Q_{p_2})$;
 \item $f_{p_3} = \one_{K_{p_3}}$;
\item $f_{x_i}$ is (essentially) a sign\footnote{See the previous footnote.} times a pseudo-coefficient of the Steinberg representation of $G(\Q_{x_i})$ for $i=1, 2, \ldots, d$;
\item Before we define the function $f^{S'}$ we explain a fact: There are only \emph{finitely many} cuspidal automorphic representations $\Pi \subset L_0^2(G^*(\Q) A(\R)^+ \backslash G^*(\A))$ of $G^*$ whose component at infinity is equal to $\Pi_\infty$ and have invariant vectors under the group $K$. In particular also the set of their possible outside $S'$-components $\Pi^{S'}$ is finite. Therefore, we find a function $f^{S'} \in \cH(G^*(\Af^{S'})) = \cH(G(\Af^{S'}))$ whose trace on $\Pi^{S'}$ is equal to $1$ if $\Pi^{S'} \cong \Pi_0^{S'}$ and whose trace equals $0$ otherwise for all $\Pi $ with $\Pi_\infty = \Pi_{0, \infty}$ and $\Pi^K\neq 0$. We fix $f^{S'}$ to be a function having this property.
\end{itemize}
The components $f_{x_i}$ are defined in the appendix~\S \ref{appendix_section2}, thus also the definition of the Hecke operator $f^{p\infty}$ is complete (see Equation~\eqref{eqn:The_Hecke_Operator}). We emphasize that at the primes $v \notin \{x_1, x_2, \ldots, x_d\}$ we take $(f_v)^{G^*(\Q_v)} = f_v$ (we have $G^*(\Q_v) = G(\Q_v)$) and at the primes $v \in \{x_1, x_2, \ldots, x_d\}$ we control the traces of the transferred function $(f_v)^{G^*(\Q_v)}$ against smooth representations via Proposition~\ref{prop:stein_coefficients}.

Due to the cuspidal component $f_{p_1}$ of $f^p$, the trace formula simplifies. Because $f_{p_2}$ is stabilizing (Labesse \cite{MR1695940}), the contribution of the proper endoscopic groups are zero, and the right hand side of Equation~\eqref{eqn:Kottwitz_Formula_Truncated} becomes a sum of the form
\begin{equation}\label{eqn:Arthur_Formula}
\sum_{\Pi} m(\Pi) \Tr( (f_\infty f^p)^{G^*(\A^p)} (\chi_b^{G(\qp)} f_\alpha), \Pi),
\end{equation}
where $\Pi$ ranges over cuspidal automorphic representations of $G^*(\A)$, and $m(\Pi)$ is the multiplicity of $\Pi$ in the discrete spectrum of $G^*(\A)$ with trivial central character on $A(\R)^+$ 
($A$ is both the split center of the group $G$ as well as the split center of the group $G^*$). Here we are applying the simple trace formula of Arthur \cite[Cor. 23.6]{MR2192011} (cf. proof of \cite[thm 7.1]{MR939691}), the correcting term in Arthur's formula vanishes due to the pseudocoefficients in the Hecke operator. The sum in Equation~\eqref{eqn:Arthur_Formula} expands to the sum
\begin{equation}\label{eqn:Automorphic_Sum_A}
\sum m(\Pi) \Tr(f_\infty^{G^*(\R)}, \Pi_\infty) \Tr(\chi_b^{G(\qp)} f_{G, \mu, \alpha}, \Pi_p) \dim \lhk (\Pi_{p_3})^{K_{p_3}} \rhk \prod_{i = 1}^d \Tr(f^{G^*(\Q_{x_i})}_{x_i}, \Pi_{x_i}), 
\end{equation}
where $\Pi$ ranges over the irreducible subspaces of $L^2_0(A(\R)^+ G^*(\Q) \backslash G^*(\A))$ such that 
\begin{itemize}
\item $\Pi^{S'} \cong \Pi_0^{S'}$; 
\item $\Pi_{p_1}$ lies in the inertial orbit $\cI(\Pi_{p_1})$ of the representation $\Pi_{p_1}$;
\item $\Pi_{p_2}$ is, up to unramified twist, isomorphic to the Steinberg representation of $G(\Q_{p_2})$;
\item $\Pi_{x_i}$ is such that $\Tr(f_{x_i}, \Pi_{x_i}) \neq 0$.
\end{itemize}
By Proposition~\ref{st:extension} we find a cuspidal automorphic representation $\pi_0$ of $G_1^*(\A)$ contained in $\Pi_0$. Let now $\Pi$ be an automorphic representation of $G^*(\A)$ contributing to Equation~\eqref{eqn:Automorphic_Sum_A}. Thus the representation $\Pi^{S'}$ is isomorphic to the representation $\Pi_0^{S'}$. Let $\pi$ be a cuspidal automorphic subrepresentation of $\Res_{[G_1^* \times Z](\A)}(\Pi)$ (Proposition~\ref{st:extension}). Enlarge $S'$ to a larger finite set $S''$ so that the representations $\pi$ and $\Pi$ are unramified for all places outside the set $S''$. At the unramified places $v \notin S''$ the representation $\Res_{[G_1^* \times Z](\Q_v)}(\Pi_{0, v})$ contains exactly one unramified representation: $\pi_{0, v}$. Therefore we have $(\pi)^{S''} \cong (\pi_0)^{S''}$. 

We now apply base change. The representation $\pi$ has the following properties:
\begin{enumerate}
\item $\pi$ is cuspidal; 
\item $\pi_\infty$ is in the discrete series;
\item $\pi_{p_1}$ is cuspidal;
\item $\pi_{p_2}$ is an unramified twist of the Steinberg representation. 
\end{enumerate}
Consider the group $G_0^{*+} := \Res_{F/F^+} (G_{0, F}^*)$. Let $\A_{F^+} := \A \otimes_\Q F^+$ and $\A_F := \A \otimes_\Q F$. Then $G_0^{*+}(\A_{F^+}) = G_0^*(\A_F)$. Because of the above properties (1), \ldots, (4), we may base change $\pi$ to an automorphic representation $BC(\pi)$ of $G_0^{*+}(\A_{F^+})$. Here we are using Corollary 5.3 from Labesse \cite{MR2856380} to see that $\pi$ has a weak base change, and then the improvement of the statement at Theorem 5.9 of [\textit{loc. cit.}], stating that\footnote{Labesse assumes that the extension $F^+/\Q$ is of degree at least $2$. We we do not have this assumption. Labesse only needs his assumption to apply the simple trace formula. For our representation $\pi$ Labesse's assumption is redundant, because we have an auxiliary place ($v = p_1$) where the representation $\pi$ is cuspidal.}, at the places where the unitary group is quasi-split (so in particular at $p$) the (local) base change of the representation $\pi_p$ is the representation $BC(\pi)_p$. By the same argument the base change $BC(\pi_0)$ exists as well. By strong multiplicity one for the group $G_0^{*+}$ we have $BC(\pi_\wp) \cong BC(\pi_{0, \wp})$ for all $F^+$-places $\wp$ above $p$.

We give the final argument when $F/F^+$ is inert at the $F^+$-place $\wp|p$, the case of the general linear groups being easier. 

The representation $\pi_{\wp}$ is of the form $\Ind_{P(F^+_\wp)}^{G_1(F_{\wp}^+)}(\rho_\wp)$ because $\pi_p$ lies in the set $\iR_1(b)$. In this induction the parabolic subgroup $P$ has Levi component $M$ with $M(\qp) = M_{\wp, 1} \times M_{\wp, 2}$ with $M_{\wp, 1}$ a general linear group and $M_{\wp, 2}$ is a unitary group. The representation $\rho_\wp$ decomposes into $\rho_\wp \cong \rho_{\wp, 1} \otimes \rho_{\wp, 2}$, where $\rho_{\wp, 1}$ is a generic unramified representation of $M_{\wp, 1}$ and $\rho_{\wp, 2}$ is an unramified twist of the Steinberg representation of $M_{\wp, 2}$. The base change is compatible with parabolic induction, the base change of a generic unramified representation is again unramified \cite{MR2856377} and the base change of a twist of the Steinberg representation is again a twist of the Steinberg representation \cite{MR2366373}. Thus the representation $BC(\pi_\wp) \cong BC(\pi_{0, \wp})$ is an induction from a representation of the form $(\chi_1, \chi_2, \ldots, \chi_{a_\wp}, \St_{\Gl_b(F^+_\wp)}, \li \chi_{a_\wp}^{-1}, \li \chi_{a_\wp -1}^{-1}, \ldots, \li \chi_1^{-1})$ where $a_\wp = \rank(M_{\wp, 1})$ and $b_{\wp} = n - a_\wp$. Consequently, $\Theta_{\pi_{0,\wp}} \circ \cN = \Theta_{\pi_{\wp}} \circ \cN$ where $\cN$ is the norm mapping from $G^{*+}_0(F_{\wp}^+)$ to $G_0^*(F_{\wp}^+)$. The norm mapping $\cN$ from $\theta$-conjugacy classes in $G^{*+}_{0}(F_{\wp}^+)$ to $G_0^*(F_\wp^+)$ is surjective for the semi-simple conjugacy classes \cite[Prop. 3.11(b)]{MR1081540}. Thus the characters $\Theta_{\pi_{\wp}}$ and $\Theta_{\pi_{0,\wp}}$ coincide on $G_0(F_\wp^+)$. By Proposition~\ref{st:atp} there is a positive constant $C_\Pi \in \R_{> 0}$ such that (for $\alpha$ sufficiently divisible) 
\begin{equation}\label{eqn:atp2}
\Tr(\chi_b^{G(\qp)} f_\alpha, \Pi_p) = C_\Pi \Tr(\chi_b^{G(\qp)} f_\alpha , \Pi_{0, p}). 
\end{equation}
Remark: To find Equation~\eqref{eqn:atp2} we applied Proposition~\eqref{st:atp} two times: first to compare $\Tr(\chi_b^{G(\qp)} f_\alpha, \Pi_p)$ with $\Tr(\chi_b^{G_1 \times Z} f_\alpha^{G_1 \times Z}, \pi_p)$, and then to compare $\Tr(\chi_b^{G(\qp)} f_\alpha, \Pi_{0,p})$ with $\Tr(\chi_b^{G_1 \times Z} f_\alpha^{G \times Z_1}, \pi_{0, p})$. 

We may now complete the proof. We return to Equation~\eqref{eqn:Automorphic_Sum_A}: 
\begin{align}\label{eqn:Automorphic_Sum_Conclusion}
\sum m(\Pi) \Tr(f_\infty^{G^*(\R)}, \Pi_{\infty}) \Tr(\chi_b^{G(\qp)} f_\alpha, \Pi_{p}) \dim \left( (\Pi_{p_3})^{K_{p_3}} \right) \prod_{i=1}^d \Tr(f^{G^*(\Q_{x_i})}, \Pi_{x_i}), 
\end{align}
where $\Pi$ ranges over the irreducible subspaces of $L_0^2(A(\R)^+ G^*(\Q) \backslash G^*(\A))$ satisfying the conditions listed below Equation~\eqref{eqn:Automorphic_Sum_A}. The following 6 facts have been established:
\begin{enumerate}
\item The sum in Equation~\eqref{eqn:Automorphic_Sum_Conclusion} is non-empty because $\Pi_0$ occurs in it (by the Propositions~\ref{st:glnonzeroclass} and \ref{st:unitnonzeroclass}, the term corresponding to $\Pi_0$ in the Sum~\eqref{eqn:Automorphic_Sum_Conclusion} is non-zero).
\item The multiplicity $m(\Pi)$ is a positive real number.
\item For any $\Pi$ in Equation~\eqref{eqn:Automorphic_Sum_Conclusion} with $\Pi \not\cong \Pi_0$ we must have $\Tr(f_\infty^{G^*(\R)}, \Pi_\infty) = \Tr(f_\infty^{G^*(\R)}, \Pi_{0, \infty})$ (here we use that $\xi$ is sufficiently regular).
\item By Equation~\eqref{eqn:atp2} the trace $\Tr(\chi_b^{G(\qp)} f_\alpha, \Pi_p)$ equals $\Tr(\chi_b^{G(\qp)} f_\alpha, \Pi_{0, p})$ up to the positive number $C_\Pi$. 
\item The dimensions $\dim \left( (\Pi_{p_3})^{K_{p_3}} \right)$ and 
$\dim \left( (\Pi_{0, p_3})^{K_{p_3}} \right)$ differ by a positive real number.
\item The product $\prod_{i = 1}^d \Tr(f^{G^*(\Q_{x_i})}, \Pi_x)$ is a non-negative real number for all automorphic representations $\Pi$ contributing to Equation~\eqref{eqn:Automorphic_Sum_Conclusion} (Proposition~\ref{prop:stein_coefficients}).
\end{enumerate}
(facts (2) and (5) are trivial). From facts (1), (2), \ldots, (6) we conclude that Equation~\eqref{eqn:Automorphic_Sum_Conclusion} must be non-zero. This completes the proof.
\end{proof}


\section{Other Shimura data}

We generalize the argument of the second chapter to more general Shimura data.

\subsection{Shimura data}\label{section_shimura_data}
Let $(G, X)$ be a Shimura datum subject to the following additional conditions:
\begin{itemize}
\item[(\textbf{C1})] We consider Shimura data in the sense of Deligne~\cite{MR0498581, MR546620}. In particular $G$ is \emph{connected}\footnote{The condition that $G$ is connected is important for our method. For non-connected groups there is not even a conjectured formula for the number of points on the Shimura variety at primes of good reduction. Furthermore non-connected groups give all sorts of technical problems in our arguments; for instance in the theory of isocrystals~\cite{MR1485921}, in the trace formula (problems with endoscopy), and also problems with the existence of (local) cuspidal representations which are used crucially at auxiliary places to simplify the trace formula. We have no idea what to expect in non-connected cases.}, and in particular this excludes PEL varieties of type $D$;
\item[(\textbf{C2})] The derived group of $G$ is simply connected;
\item[(\textbf{C3})] The maximal $\Q$-split torus in the center of $G$ coincides with the maximal $\R$-split torus in the center of $G$.
\end{itemize}

For Shimura varieties satisfying conditions (C1), (C2) and (C3) we prove, assuming a condition on the group, that the formula of Kottwitz does not vanish if one restricts it to the contribution of a certain Newton stratum. \emph{Conjecturally} this implies that the Newton strata are non-empty. More precisely, non-emptiness follows if one assumes the following two conjectures:
\begin{itemize}
\item[(\textbf{C4})] The Rapoport-Langlands conjecture\footnote{Whenever we mention this conjecture, we mean the version as refined by Milne (See \S~\ref{sect:raplan}).} \footnote{In fact, we do not the full conjecture. We need equality of Equations~\eqref{pointformula_raplankottwitz} and~\eqref{pointformula_integrals} below, which is a weaker statement.} is true for $(G, X)$;
\item[(\textbf{C5})] A technical~---and unfortunately deep---~ hypothesis (Hyp.~\ref{hypG}) on the cuspidal automorphic spectrum of the quasi-split inner form $G^*$ of $G$.
\end{itemize}
In the next chapters 4 and 5 we prove condition (C5) for some groups, so that in those cases our result is only conditional on (C4). More precisely, we prove:
\begin{itemize}
\item In Chapter 4 we show using Arthur's results~\cite{arthurbook} that (C5) is satisfied for symplectic groups of similitudes. For Shimura varieties of PEL type (C), the results of Kottwitz in~\cite{MR1124982} then prove enough of (C4) to conclude that the Newton strata are non-empty. 
\item In Chapter 5 we show using~\cite{arthurbook} that (C5) is essentially\footnote{See the third remark below Theorem~\ref{theorem_gspin} for the precise statement.} true for certain spinor groups of type B. This provides us with sufficient information to conclude that the Newton strata are nonempty, conditional on (C4). 
\end{itemize}
We should mention that the results in~\cite{arthurbook} are currently still conditional on the stabilization of the twisted trace formula for $\GL_n$. Thus our result for PEL varieties of type C and the odd spinor groups is conditional on this stabilization as well. 

To clarify our assumptions: We assume throughout this chapter that $(G, X)$ satisfies (C1), (C2) and (C3), but not necessarily (C4) and (C5). 

\subsection{Notations/Conventions}\label{situation_good_reduction}
Let $(G, X)$ be a Shimura datum satisfying (C1)..(C3). We fix the following notations and conventions:
{\small
\begin{itemize}
\item $E$ is the reflex field of $(G, X)$;
\item $p$ is a prime number such that the group $G$ and the field $E$ are unramified at $p$;
\item $\p$ is a (any) prime of $E$ above $p$;
\item $\fq$ is the residue field of $E$ at $\p$;
\item $E_{\p}$ is the completion of $E$ at $\p$; 
\item $\alpha$ is a positive integer;
\item $E_{\p, \alpha}$ is the unramified extension of degree $\alpha$ of $E_{\p}$;
\item $\fqa$ is the residue field of $E_{\p, \alpha}$;
\item $L$ is the completion of the maximal unramified extension of $\qp$, and we fix embeddings $E_{\p, \alpha} \subset L$;
\item $K_p \subset G(\qp)$ is a hyperspecial subgroup;
\item $K^p \subset G(\afp)$ is sufficiently small so that the stack $\Sh_K/E$ is smooth and represented by a scheme;
\item $\Sh_K/\cO_E$ is the (conjectural) canonical model of Milne~\cite{MR1155229}; its existence is proved by Vasiu and Kisin~\cite{MR2669706, MR2541427} for Shimura data of Abelian type;
\item we assume $K^p$ is sufficiently small so that the canonical model $\Sh_K$ is smooth. 
\end{itemize}
}

\subsection{The Rapoport-Langlands conjecture}\label{sect:raplan}
Briefly, the Rapoport-Langlands conjecture\footnote{Milne made important refinements~\cite{MR1155229} to the original conjecture of Rapoport-Langlands~\cite{MR895287}. In the original work of Rapoport-Langlands it is not explained how the variety $S/E$ is to be reduced modulo primes of good reduction. Rapoport and Langlands simply state that there should exist a model for which Equation~\eqref{raplanmilne} is the set of points over finite fields. Milne refined the conjecture by first conjecturing the existence of a canonical model, and then conjecturing that for this canonical model, Equation~\eqref{raplanmilne} is the set of points. Milne also extends the conjecture to groups whose derived group is not simply connected.} (cf. \cite{MR1155229, MR895287}) states
\begin{equation}\label{raplanmilne}
S(\lfq) = \coprod_{\varphi} I(\varphi) \backslash X(\varphi) /K^p G(\cO_L), 
\end{equation}
where
\begin{itemize}
\item $\varphi$ ranges over the special morphisms of the Tannakian category $\Rep_\Q(G)$ to the category of $\fqa$-motives $\Mot(\fqa)$. 
\item Each such morphism $\varphi$ is a motive with additional $G$-structure and $I(\varphi)$ is the automorphism group of $\varphi$.
\item $X(\varphi)$ is a set of lattices in the product over all the $\ell$-adic realizations ($\ell \neq p$) of $\varphi$, and lattices in the crystalline realization of $\varphi$.
\end{itemize}

\subsection{Truncated Kottwitz formula}\label{sect:Truncated_Kottwitz_Formula}
Assume that the Rapoport-Langlands-Milne conjecture is true for our Shimura datum $(G, X)$. Let $x\in S(\fqa)$ be a point and let $\varphi$ be the parameter corresponding to $x$ in Equation~\eqref{raplanmilne}. Write $\varphi_p$ for the crystalline realization of $\varphi$. Fix an isocrystal $b \in B(G_{\qp}, \mu)$, then we can define
$$
S^b(\lfq) := \coprod_{\bo{\quad\varphi}{ b \cong \varphi_p \otimes L}} I(\varphi) \backslash X(\varphi) /K^p G(\cO_L), 
$$
so we take the disjoint union ranging over the same morphisms as in Equation~\eqref{raplanmilne}, but we restrict to those whose crystalline realisation is, up to isogeny, equal to $b$. We take this as definition for $S^b(\lfq)$. The subset $S^b(\fqa) \subset S^b(\lfq)$ is the set of elements $x$ that are invariant under $\Phi_\p^\alpha$.

\begin{remark}
We defined $S^b(\fqa)$ \emph{only as a set}; we did not give $S^b$ the structure of a scheme and we did not prove that the varieties $S^b \subset S$ are locally closed and form a stratification of $S$. Vasiu shows in~\cite[Thm.~5.3.1]{MR2789744} that the Newton strata form locally closed varieties for Shimura varieties of Hodge type.
\end{remark}

The decomposition in Equation~\eqref{raplanmilne} is that of a set with additional operators~\cite{MR1155229}. Let $g \in G(\Afp)$ and put $K_g^p := K^p \cap g K^p g^{-1}$. Let $a \colon S_{K^p_g} \to S_K^p$ be the mapping induced by $g$ and let $c$ be the usual covering map of $S_{K_g^p}$ onto $S_{K^p}$. Then the couple $(a, c)$ induces a correspondence $f^{p\infty}$ on the Shimura variety $S_K$. We can restrict the correspondence $f^{p\infty}$ to the set of geometric points of the different Newton strata $S^b(\lfq)$, $b \in B(G_{\qp}, \mu)$. 

Let $\xi$ be an irreducible representation of the group $G_\C$. Let $\ell \neq p$ be a prime number. Fix an isomorphism $\iota$ between $\C$ and $\lql$. Then $(\iota, \xi)$ induce an $\ell$-adic local system $\cL$ on $S$. Let $\Fix(f^{p\infty} \times \Phi_{\p}^\alpha)$ be the set of fix points of $f^{p\infty} \times \Phi_{\p}^{\alpha}$ acting on the scheme $S_K$. The set of fix points of $f^{p\infty} \times \Phi_{\p}^\alpha$ acting on $S^b(\lfq)$ is the set of $x' \in \Fix(f^{p\infty} \times \Phi_{\p}^\alpha)(\lfq)$ such that the image $x$ of $x'$ in $S(\lfq)$ lies in the subset $S^b(\lfq) \subset S(\lfq)$. We write $\Fix^b(f^{p\infty} \times \Phi_{\p}^\alpha)(\lfq)$ for this set. We emphasize egain that we did not define $\Fix^b(f^{p\infty} \times \Phi_{\p}^\alpha)$ as a scheme; we only defined the set of geometric points. We consider the expression\footnote{Using the isomorphism $\iota \colon \C \isomto \lql$ we view the number in Equation~\eqref{pointformula_raplankottwitz} as a complex number. More precisely, the function $f^{p\infty}$ is made $\lql$-valued via $\iota$, then we compute the trace of $f^{p\infty} \times \Phi_\p^\alpha$ on $\cL_x$, giving a number in $\lql$ which in turn we send to $\C$ via the isomorphism $\iota^{-1}$.}
\begin{equation}\label{pointformula_raplankottwitz}
\sum_{x' \in \Fix^b(f^{p\infty} \times \Phi_\p^\alpha)(\lfq)} \Tr(f^{p\infty} \times \Phi_\p^\alpha, \cL_x), 
\end{equation} 
where the point $x \in S_K(\lfq)$ is the image of $x' \in \Fix^b(f^{p\infty} \times \Phi_\p^\alpha)(\lfq)$ in the Shimura variety. The argument of Milne shows~\cite{MR1155229}, using Milne's formulation of the Rapoport-Langlands conjecture, that Formula~\eqref{pointformula_raplankottwitz} equals the familiar orbital integral expression
\begin{equation}\label{pointformula_integrals}
\sum_{(\gamma_0; \gamma, \delta)} c(\gamma_0; \gamma, \delta) O_\gamma(f^{p\infty}) \TO_\delta(\chi_b^G \phi_\alpha) \Tr \xi(\gamma_0), 
\end{equation} 
with $(\gamma_0; \gamma, \delta)$ ranging over the Kottwitz triplets. The derived group of $G$ is simply connected, and because $F^+$ is totally real, the maximal $\R$-split torus in $G_\R$ coincides with the maximal $\Q$-split torus in $G$. Thus the conditions in Kottwitz's stabilization argument are satisfied (cf. \cite[Sect. 1]{MR1124982}). Therefore, Expression~\eqref{pointformula_integrals} is equal to
\begin{equation}\label{equnow34}
\sum_{H \in \cE} \iota(G, H) \ST_e^H( h^{p\infty} (\chi_b^G \phi_\alpha)^H),  
\end{equation}
where the functions $h^{p\infty} \in \cH(H(\Afp))$ are defined by Kottwitz~\cite{MR1124982}. 

Now \emph{forget} the assumption that the Rapoport-Langlands-Milne conjecture is true for $(G, X)$. Formula~\eqref{pointformula_integrals} still makes sense, and is still equal to Formula~\eqref{equnow34}. The goal of this chapter is to show that Formula~\eqref{equnow34} is non-zero for particular choices of $(f^{p\infty}, \cL)$. Consequently Formula~\eqref{pointformula_integrals} is non-zero as well for those $(f^{p\infty}, \cL)$. \emph{Conditional} on the Rapoport-Langlands conjecture, this implies that Formula~\eqref{pointformula_raplankottwitz} is non-zero, and thus that the $b$-th Newton stratum is non-empty.

\renewcommand{\Pr}{\textup{Pr}}

\subsection{Local computations}
We construct a certain class $\iR(b)$ of representations of local reductive groups associated to admissible isocrystals $b \in B(G_{\qp}, \mu)$. 

We drop the global notations that have been introduced earlier and revert to the local setting. Thus, let $G$ be a quasi-split connected reductive group over $\qp$. We reuse the notations $P_0$, $T$, $Z$, $A$, $Z$, $\ldots$ introduced in Chapter 1. On top of those notations, we introduce several more: 
{\small \begin{itemize}
\item $E \subset \lqp$ is a finite unramified extension of $\qp$.
\item for each positive integer $\alpha$, $E_\alpha \subset \lqp$ is the unramified extension of $E$ of degree $\alpha$. 
\item $\mu \in X_*(T)$ is a miniscule dominant cocharacter. 
\end{itemize} }
Below we compute both with the relative root system of $G$, and the relative root system of the group $G_{E_{\alpha}}$. To simplify, we assume $\alpha$ is sufficiently divisible so that $G_{E_{\alpha}}$ splits. Then the relative root system of $G_{E_{\alpha}}$ is the absolute root system of $G$. We need the following notations:
{\small \begin{itemize}
\item Let $P \subset G$ a parabolic subgroup. We write $A_P^+$ for the $E_\alpha$-split center of the group $P_{E_\alpha}$. 
\item $A_0^+ := A_{P_0}^+ = T_{E_\alpha}$.
\item $\iA_P := X_*(A_P^+)_\R$. 
\item If $P \subset Q$ is an inclusion of standard parabolic subgroups, then $\iA_P^Q:= \iA_P/\iA_Q$. 
\item $W_G^+$ is the Weyl group of $A_{P_0}^+$ in $G_{E_{\alpha}}$.
\item $\Delta_0^+$ is the set of simple roots on $\iA_0 := \iA_{P_{0, E_\alpha}}$ corresponding to the Borel subgroup $P_{0, E_\alpha}$. 
\item To the parabolic subgroup $P$ we associate the subset $\Delta_P^+ \subset \Delta^+$ consisting of those roots acting non trivially on $A_P$. Define $\ia_P$ to be $X_*(A_P)_\R$,
\item The map $\pi_P$ is the projection $\iA_0 \surjects \iA_P$.
\item The map $\Pr_P$ is the projection $\iA_0 \surjects \iA^P_0$. 
\item If $X \in \iA_0$, then $X = \pi_P(X) + \Pr_P(X)$.
\end{itemize} }
We make a convention on the center and the fundamental weights: The $\varpi_\alpha$ are defined as the basis of $(\ia_0^G)^*$ dual to the basis of coroots. In particular, the $\varpi_\alpha$ do not make sense on the space $\ia_0$ (they only make sense \emph{after} we quotient $\ia_0$ out by the center $\ia_G$). We extend the definition of the $\varpi_\alpha$ to the whole space $\ia_0$, by defining $\varpi_\alpha(X) = 0$ for all $X \in \ia_G$. We extend in the same way the fundamental weights $\varpi_\alpha$, $\alpha\in \Delta_0^+$ from $\iA_0^G$ to $\iA_0$. 

We can project downwards from the absolute root system $\iA_0$ to the relative root system $\ia_0$:
{\small \begin{itemize}
\item $\theta \in \Gal(E_\alpha/\qp)$ is a generator of the Galois group. 
\item For each $P$ we identify $\ia_P$ with the set of $\theta$-invariant elements in $\iA_P$. 
\item The norm mapping $\cN$ from $X_*(T)$ to $X_*(A_0)$ is given by $\cN(X) = \sum_{\sigma \in \Gal(E_{\alpha}/\qp)} \sigma(X)$.
\item The maps $\cN$ and $\pi_P$ commute: For all $X \in \iA_0$ we have $\cN(\pi_P(X)) = \pi_P(\cN(X))$.
\item The maps $\cN$ and $\Pr_P$ commute as well.
\end{itemize}}

We start with a lemma on the base change morphism

\newcommand{\Frob}{\textup{Frob}}

\begin{lemma}
Let $\BC \colon \cH(G(E_{\alpha})) \to \cH(G(\qp))$ be the base change mapping. Then $\BC$ corresponds via the Satake isomorphisms to the norm mapping:
\begin{align*}
\cN \colon \C[X_*(A_{0, \alpha})]^{W_\alpha} \rightarrow \C[X_*(A_0)]^W, \ X_*(A_{0, \alpha}) \owns X \mapsto \sum_{\sigma \in \Gal(E_{\alpha}/\qp)} \sigma(X).
\end{align*}
\end{lemma}
\begin{proof}
The base change mapping $\BC \colon \cH_0(\Res_{E_\alpha/\qp} G_{E_\alpha}) \rightarrow \cH_0(G)$ is compatible with the constant term morphism $f \mapsto f^{(P_0)}$. Therefore the lemma reduces to the case where $G$ is commutative. If $G$ is commutative, then base change of representations is simply the composition of smooth characters $\psi$ of $G$ with the norm mapping $\cN$ of $E_\alpha/\qp$. Let $X \in X_*(G_{E_\alpha})$, then $X$ is a cocharacter, so an algebraic morphism $\qpm\rightarrow G(E_\alpha)$, which we can evaluate at $p$. By the definition of the base change map we have $(\psi \circ \cN)(X(p)) = \psi(\BC(X)(p))$ for all characters $\psi$ of $G$. We compute
$$
(\psi \circ \cN)(X(p)) = \psi (( \sum_{\sigma \in \Gal(E_\alpha/\qp)} \sigma(X)) (p)). 
$$
Since this equation holds for all $\psi$, we get $\BC(X) = \sum_{\sigma \in \Gal(E_\alpha/\qp)} \sigma(X)$. This proves the lemma.
\end{proof}

Let $b$ be a $\mu$-admissible isocrystal with $G$-structure. Write $\nu := \li \nu_b$ for the slope morphism of $b$. We associate to $b$ a certain class of representations $\iR(b)$ of the group $G$:

\begin{definition}\label{therbdef}
Let $P = MN$ is the standard parabolic subgroup of $G$ contracted by $b$. Thus $\Delta_P \subset \Delta_0$ is the subset of $\alpha \in \Delta_0$ such that $\langle \alpha, \nu \rangle > 0$. We define another standard parabolic subgroup $Q = LU \subset P$. Let $\Delta_Q^P$ be the set of $\alpha \in \Delta_0^P$ such that there exists $w \in W_G^+$ such that $\langle \varpi_\alpha , \cN(w(\mu) - \nu) \rangle = 0$. Then $Q \subset P$ is the standard parabolic corresponding to $\Delta_Q^P$. Let $\iR(b)$ be the set of $G$-representations $\Pi$ of the form $\Ind_Q^G[\St_L(\chi)]$, where $\St_L$ is the Steinberg representation of $L$, and $\chi$ ranges over the unramified unitary characters of the group $L$ such that the following condition holds. For certain $\chi$ the induction $\Ind_Q^G[\St_L(\chi)]$ may be reducible, we keep only those $\chi$ such that $\Ind_Q^G[\St_L(\chi)]$ is irreducible. 
\end{definition}

\begin{example}
To help the reader get a feeling for the above definition, let us explain what it gives for a concrete example. Assume (for the moment) that $G$ is the symplectic group of similitudes over $\Q$, and $\mu$ is the usual $(1_g, 0_g)$ cocharacter coming from the classical Shimura varieties $\cA_g$ classifying polarized Abelian varieties of dimension $g$. Giving an isocrystal in $B(G_{\qp}, \mu)$ is the same as giving a sequence of slopes $\lambda_1 \geq \lambda_2 \geq \ldots \geq \lambda_{2g} \in \Q$ with, for all $i$, $0 \leq \lambda_i \leq 1$, $\lambda_{2g+1-i} + \lambda_i = 1$ and such that the Newton polygon $\cG_\lambda$ associated to $\{\lambda_i\}$ only has breakpoints at points $(x,y) \in \Q^2$ for which both coordinates $x,y$ are integers. It is important to note that, even though all the breakpoints have integral coordinates, there could be some integral points on $\cG_\lambda$ which are not breakpoints. Mark all the integral points (including the breakpoints) on the Newton polygon $\cG_\lambda$. These markings define a composition of $n = 2g$, and in turn this composition defines a standard parabolic subgroup $Q$ of $G_{\qp}$. This $Q$ is the same as the $Q$ which occurred in  Definition~\ref{therbdef}.  
For other split classical groups whose derived group is simply connected the idea remains the same: The $\mu$-admissible isocrystals are classified by sequences of slopes $\{\lambda_i\}$ satisfying a certain symmetry, and to get the parabolic subgroup $Q \subset G_\qp$ you mark all the integral points on the Newton polygon $\cG_{\lambda}$ associated to $\{\lambda_i\}$. 
\end{example}

\begin{remark}
Assume $G$ is one of the unitary groups of similitudes from the second chapter. Then, the class of representations $\iR(b)$ that we consider here is slightly different from the class of representations that we considered in Chapter 2. (So we claim that the isolation argument works for both classes of representations). 
\end{remark}

\begin{proposition}
Let $\Pi \in \iR(b)$, then for $\alpha$ sufficiently divisible, we have $\Tr(\chi_b^G f_\alpha, \Pi) \neq 0$.
\end{proposition}
\begin{proof}
To avoid confusion we change here (and only here) the symbol $\alpha$ in $f_\alpha$ to a $\beta$, so $f_\beta$ is the degree $\beta$ Kottwitz function.

Let $a_X \in \C$ be the coefficients of the Kottwitz function $f_\beta$ occurring in its Satake transform: $$\cS(f_\beta) = q^{-\beta \langle \rho_G,\mu\rangle} \sum_{X \in X_*(T)} a_X \cdot X \in \C[X_*(T)]^{W}.$$ For all $X$ we have $a_X \in \Z_{\geq 0}$. We have $\Tr(\chi_b^G f_\beta, \pi) = \Tr(\chi_b^G f_\beta^{(Q)}, \St_L(\chi))$. We have $\chi_b^G|_L = (\chi_c^M \cdot \eta_b)|_L$ by Lemma~\ref{st:basicstruncatedcompact}. By applying Proposition~\ref{st:truncatedtraceonsteinberg} to the group $L$ we get
\begin{equation}\label{dimension.CC}
\Tr((\chi_c^M \eta_b)|_L f_\beta^{(Q)}, \St_L(\chi)) = \left. \sum_X a_X \cdot X \right|_{\chi \delta_{P_0 \cap L}^{-1/2}}, 
\end{equation}
where $X \in X_*(A_0)$ ranges over the monomials such that 
\begin{itemize}
\item[(\textbf{S1})] $\pi_Q(X) = \beta \nu$; 
\item[(\textbf{S2})] $\langle \Pr_Q(X), \varpi_\alpha \rangle > 0$ for all $\alpha \in \Delta_0^Q$.
\end{itemize}
To show that Equation~\eqref{dimension.CC} is non-zero we may forget about the twist $\chi$, and it suffices to show that the evaluation $\left. \sum_X a_X \cdot ( X \right|_{\delta_{P_0 \cap L}^{-1/2}} ) $ is non-zero. Because for all $X$, $a_X \geq 0$, it suffices to show that the sum is non-empty, \ie, that there exists at least one monomial $X$ satisfying (S1), (S2) and $a_X \neq 0$. By definition of the Kottwitz function we have $\cS(\phi_\alpha) = q^{\alpha \langle \rho_G, \mu \rangle} \sum_{X \in W_G^+ \cdot \mu} [X]$. By a theorem of Chai~\cite[Thm.~5.1]{MR1781927} there exists an $X$ in the Weyl group orbit $W_G^+\cdot\mu$ such that the projection of $X$ to the chamber $\iA_P$ coincides with $\nu$. We consider this $X$, and we show that its norm $Y:=\cN(X)$ satisfies the conditions (S1) and (S2). Condition (S1): We have $\pi_P(X) = \nu \in \iA_P$. Apply the norm $\cN$ to get $\cN(\pi_P(X)) = \cN(\nu) = \beta \nu$. Since $\cN(\pi_P(X)) = \pi_P(\cN(X))$, we get $\pi_P(Y) = \beta \nu$  (observe that $\theta(\nu) = \nu$ for all $\theta \in \Gal(E_\beta/\qp)$, so $\cN(\nu) = \beta \nu$). Projecting to $Q$ we get $\pi_Q (Y) = \pi_Q\pi_P(Y) = \pi_Q (\beta\nu) = \beta\nu$. For the equality $\pi_Q(\nu) = \nu$, we used that $P$ is contracted by $\nu$ and therefore $\nu$ lies in the center of $\iA_P$. Thus (S1) is satisfied. We now show (S2): Because $Q \subset P$, we get from $\pi_P(X) = \nu$ that $\pi_Q(X) = \pi_Q(\pi_P(X)) = \pi_Q(\nu) = \nu$. Write $X = \Pr_Q(X) + \pi_Q(X) = \Pr_Q(X) + \nu$. For all $\alpha \in \Delta_0^+$, we have $\langle \varpi_\alpha, X \rangle \geq \langle \varpi_\alpha, \pi_P(X) \rangle$ \cite[3.6]{MR1781927}, and thus $\langle \varpi_\alpha, X \rangle \geq \langle \varpi_\alpha, \nu\rangle$, and $\langle \varpi_\alpha, X - \nu \rangle \geq 0$. For all $\theta \in \Gal(E_\alpha/\qp)$ we have $\theta(\varpi_\alpha) = \varpi_{\theta \alpha}$ and thus $0 \leq \langle \theta(\varpi_\alpha), X - \nu \rangle = \langle \varpi_\alpha, \theta(X) - \nu \rangle$. By summing over all $\theta$, we get $\langle \varpi_\alpha, Y - \beta \nu \rangle \geq 0$. By definition of $Q$, the inequality $\langle \varpi_\alpha, Y - \beta \nu \rangle \geq 0$ is \emph{strict} if $\alpha \in \Delta_0^L$. Because $X - \nu = \Pr_Q(X)$, (S2) is true as well.
\end{proof}

\subsection{A technical Hypothesis on $G$}
We make a certain technical assumption on the group $G$ for which we account in later chapters in case $G$ is a symplectic group  of similitudes or $G$ is an odd spinor group.

Let $p$ be a prime number.

\begin{hypothesis}\label{Ghyp}\label{hypG}
Let $G^*$ be the quasi-split inner form of the group $G$. Let $\Pi, \Pi_0$ be two cuspidal automorphic representations of the group $G^*(\A)$ satisfying the following conditions:
\begin{itemize}
\item $\Pi_{0,\infty}$ is $\xi$-cohomological for a sufficiently regular algebraic representation $\xi$ of $G_\C$ whose central character is trivial; 
\item $\Pi_\infty \cong \Pi_{0, \infty}$;
\item Let $S$ be the set of prime numbers $x$ where $G_{\Q_x}$ is ramified. Then, for each $x \in S$, the representation $\Pi_x, \Pi_{0, x}$ are unramified twists of the Steinberg representation of $G^*(\Q_x)$. 
\item Let $p_1, p_2$ be two prime numbers different from $p$, $p_i \notin S$ and $G_{\Q_{p_i}}$ is split for $i=1,2$. The representations $\Pi_{p_1}, \Pi_{0, p_1}$ are unramified twists of the Steinberg representation. The representations $\Pi_{p_2}, \Pi_{0, p_2}$ are cuspidal, and unramified twists of each other.
\item For all primes $x \notin S \cup \{p, p_1, p_2\}$ we have $\Pi_x \cong \Pi_{0, x}$. 
\end{itemize}
Then the assumption we make on $G^*$ is the following implication:
$$
\Pi_{0, p} \in \iR(b) \Longrightarrow \Pi_p \cong \Pi_{0, p}(\chi), 
$$
for some unramified character $\chi$ of $G(\qp)$ of finite order.
\end{hypothesis}

\subsection{Back to the isolation argument}
We show non-emptiness of the Newton strata\footnote{More precisely, we show non-vanishing of Equation~\eqref{equnow34}, which, conditional on the Rapoport-Langlands conjecture, or a suitable form of the Kottwitz formula, implies non-emptiness of the Newton strata.}, \emph{conditionally} on our ugly hypothesis on $G$. 

\begin{theorem}\label{conditional_theorem}
Let $(G, X)$ be a Shimura datum satisfying conditions (C1), (C2) and (C3) from Section~\ref{section_shimura_data}. Let $p$ be a prime number where the Shimura variety associated to $(G, X)$ has good reduction (precise conditions: Section~\ref{situation_good_reduction}). Assume Hypothesis~\ref{Ghyp} is true for $(G, p)$. Let $b$ be a $\mu$-admissible $G_{\qp}$-isocrystal. Then Formula~\eqref{equnow34} is non-zero for a suitable choice of operator $f^p$ and local system $\cL$.
\end{theorem}
\begin{proof}
Let $\alpha$ be any positive integer. Pick a prime $\ell \neq p$ and fix an isomorphism $\lql \cong \C$ (and suppress it from all notations). Let $\xi$ be an irreducible complex (algebraic) representation of $G$, and $\cL$ is the corresponding $\ell$-adic local system on the Shimura tower. We start with Formula~\eqref{equnow34} from Section~\ref{sect:Truncated_Kottwitz_Formula}:
\begin{equation}\label{eqn:Symp_Kottwitz_restricted} 
\sum_{H \in \cE} \iota(G, H) \ST_e((\chi_b^G f_\alpha)^H)
\end{equation}
where $\cE$ is the (finite) set of endoscopic groups $H$ associated to $G$ and unramified at all places where the data $(G, K)$ are unramified. We show that this formula is non-zero for certain $f^p, \cL$ via a variant on the proof of Theorem~\ref{st:mainthmA}. Pick auxiliary prime numbers $p_1, p_2$ where the group $G(\Q_{p_i})$ is split. Let $\Pi_{p_2}^\sharp$ be a cuspidal representation of the group $G(\Q_{p_2})$ \cite{kretexistence}. Let $S$ be the set of primes $x$ where the group $G(\Q_x)$ ramifies. We take the compact open group $K \subset G(\Af)$ with the following properties:
\begin{itemize}
\item $K$ is a product $\prod_v K_v \subset G(\Af)$ of compact open groups;
\item $K$ is sufficiently small so that $\Sh_{K}$ is smooth;
\item for all $x \notin S \cup \{p_1, p_2\}$ the group $K_x$ is hyperspecial;
\item the group $K_p$ is an Iwahori subgroup of $G(\qp)$;
\item the group $K_{p_1}$ is an Iwahori subgroup of $G(\Q_{p_1})$;
\item the group $K_{p_2}$ is sufficiently small so that $(\Pi_{p_2}^\sharp)^{K_{p_2}} \neq 0$;
\item for each $x \in S$, the group $K_x$ is sufficiently small so that Steinberg representation of $G(\Q_x)$ has vectors invariant under $K_x$; 
\end{itemize}
The group $K^* \subset G^*(\Af)$ is a (any) group with the following properties:
\begin{itemize}
\item $K^*$ is a product $\prod_v K_v^* \subset G^*(\Af)$ of compact open groups;
\item for any prime $v \notin \{x_1, \ldots, x_d\}$ we have $K_v^* = K_v \subset G(\Q_v) = G^*(\Q_v)$;
\item for all $x \in S$, the Steinberg representation of $G^*(\Q_x)$ has vectors that are invariant under $K_x$.
\end{itemize}
We claim that there exists a cuspidal automorphic representation $\Pi_0$ of $G^*(\A)$ such that
\begin{itemize}
\item $\Pi_{0, \infty}$ is $\xi$-cohomological for an irreducible representation $\xi$ of $G_\C$ which is sufficiently regular. 
\item $\Pi_0^{K^*} \neq 0$; 
\item $\Pi_{0, p} \in \iR(b)$;
\item $\Pi_{0, p_1}$ is an unramified twist of the Steinberg representation of $G(\Q_{p_1})$;
\item $\Pi_{0, p_2}$ lies in the inertial orbit of the representation $\Pi_{p_2}^\sharp$;
\item $\Pi_{0, S}$ is a twist of the Steinberg representation of $G^*(\Q_S)$. 
\end{itemize}
We will prove the claim at the end of this section. 

We consider functions $f^{\infty p} = \bigotimes_v f_v \in \cH(G(\Af))$ where
\begin{itemize}
\item $f_{p_1}$ is a sign (cf. Prop.~\ref{prop:stein_coefficients}) times a pseudo-coefficient of the Steinberg representation of $G(\Q_{p_1})$;
\item $f_{p_2}$ is a pseudo-coefficient on $G(\Q_{p_2})$ of the representation $\Pi_{p_2}$;
\item for each $x \in S$, $f_x$ is a sign (cf. Prop.~\ref{prop:stein_coefficients}) times a pseudo-coefficient of the Steinberg representation of $G(\Q_x)$;
\item $f^{S \cup \{p, p_1, p_2\}} \in \cH(G^*(\Af^{S \cup \{p, p_1, p_2\}}))$ is a function such that, for all $\Pi$ occurring in $L_0^2(G^*(\Q) A(\R)^+ \backslash G^*(\A))$ with $\Pi^K \neq 0$ and $\Pi_\infty$ is $\xi$-cohomological, we have 
$$
\Tr(f^{S \cup \{p, p_1, p_2\}}, \Pi^{S \cup \{p, p_1, p_2\}}) = 1,
$$ if $\Pi^{S \cup \{p, p_1, p_2\}} \cong \Pi_0^{S \cup \{p, p_1, p_2\}}$ and $\Tr(f^{S \cup \{p, p_1, p_2\}}, \Pi^{S \cup \{p, p_1, p_2\}}) = 0$ otherwise. 
\end{itemize}
For the component $f_\infty$ of $f$ at infinity we also make a choice: Let $P_\infty$ be the the packet consisting of those discrete series representations $\Pi_\infty$ of $G^*(\R)$ whose infinitesimal character and central character coincide with those of $\xi^\vee$. We have $\Pi_{0, \infty} \in P_\infty$. The stable orbital integrals of the function $f_\infty$ are prescribed by the stabilization argument of Kottwitz~\cite[Eq.~(7.4)]{MR1044820}. Kottwitz has to show in [\textit{loc. cit}] that functions with these stable integrals exist. Let for each $\Pi_\infty \in P_\infty$, $f_{\Pi_\infty}$ be a Clozel-Delorme psuedocoefficient for $\Pi_\infty$. If $H = G^*$ is the maximal endoscopic group of $G$, then Kottwitz takes $f_\infty$ to be the average of the functions $f_{\Pi_\infty}$ with $\Pi_\infty$ ranging over $P_\infty$. The stable orbital integrals of the functions $f_{\Pi_\infty}$ are all the same (we explain this below). Therefore, on the maximal endoscopic group, one may also take $f_\infty = f_{\Pi_{0, \infty}}$. Since the stable orbital integrals of $f_\infty$ and $f_{\Pi_{0, \infty}}$ coincide, the stabilization argument of Kottwitz still carries through when the role of $f_\infty$ is replaced by that of $f_{\Pi_{0, \infty}}$. In the rest of the argument we put $f_\infty := f_{\Pi_{0, \infty}}$, so that $\Tr(f_\infty, \Pi_\infty) \neq 0$ implies $\Pi_\infty \cong \Pi_{0, \infty}$ for any automorphic representation $\Pi$ in the discrete spectrum of $G^*(\R)$.

To see that the stable orbital integrals of $f_{\Pi_\infty}$ do not depend on $\Pi_\infty \in P_\infty$, one can appeal to a result of Kottwitz~\cite[lemma 3.6]{MR1163241}\footnote{This article is written for Shimura varieties associated to division algebras, but the lemma we are referring to is true in general (for connected real reductive groups).} saying that for any stable elliptic conjugacy class $\gamma_\infty$ in $G(\R)$, $\SO_{\gamma_\infty}(f_{\Pi_\infty})$ is equal to $\Tr \xi(\gamma_\infty) \cdot \textup{vol}(I) \cdot e(I)$, where $\textup{vol}(I)$ is some volume term, and $e(I)$ is a sign associated to the inner form $I$ of the centralizer of $\gamma_\infty$ in $G^*(\R)$ that is anisotropic modulo the center. The stable orbital integral $\SO_{\gamma_\infty}(f_{\Pi_\infty})$ vanishes in case $\gamma_\infty$ is not elliptic. Thus the expression for $\SO_{\gamma_\infty}(f_{\Pi_\infty})$ does not depend on $\Pi_\infty$ (it only depends on $\gamma_\infty$ and the representation $\xi$).  

For this choice of Hecke function $f$, the right hand side of Equation~\eqref{eqn:Symp_Kottwitz_restricted} simplifies to 
\begin{equation}\label{eqn:Symp_expanded}
\sum_\Pi m(\Pi) \Tr(f_\infty, \Pi_\infty) \Tr(\chi_b^{G(\qp)} f_{G\mu\alpha}, \Pi_p) \Tr(f_{S \cup \{p_1\}}, \Pi_{S \cup \{p_1\}})
\end{equation}
where $\Pi$ ranges over the irreducible subspaces of $L_0^2(Z(\R)^+ G(\Q) \backslash G(\A))$. Let $\Pi$ be an automorphic representation of $G^*$ contributing to Formula~\eqref{eqn:Symp_expanded}. We show that the trace of $f$ on $\Pi$ is of the same sign as the trace of $f$ on $\Pi_0$. By Prop.~\ref{prop:stein_coefficients} we have $\Tr(f_{S \cup \{p_1\}}, \Pi_{S \cup \{p_1\}}) \in \R_{\geq 0}$. Due to our choice of $f$ we know that
\begin{itemize}
\item $\Pi_{p_i} \cong \Pi_{0, p_i}(\chi_{p_i})$, where $\chi_{p_i}$ is an unramified character of finite order ($i = 1,2$);
\item $\Pi_x \cong \Pi_{0, x}(\chi_x)$ where $\chi_x$ is a character of finite order of $G(\Q_x)$, and $x \in S$;
\item $\Pi_\infty \cong \Pi_{0, \infty}$; 
\item $\Pi_\ff^{S \cup \{p, p_1, p_2\}} \cong \Pi_{0, \ff}^{S \cup \{p, p_1, p_2\}}$.
\end{itemize}
(The characters $\chi_{p_i}$ and $\chi_x$ are of finite order because $\xi$ has trivial central character.)

In particular $\Pi$ is a cuspidal automorphic representation. By construction, the couple ($\Pi$, $\Pi_0$) satisfies the conditions in Hypothesis~\ref{Ghyp}. Thus, by the conclusion of this Hypothesis, the representation $\Pi_p$ is a finite unramified twist of $\Pi_{0, p}$. We return to Equation~\eqref{eqn:Symp_expanded}. Because there are only finitely many contributing $\Pi$, we can take $\alpha$ sufficiently divisible so that $\Tr(\chi_b^{G(\qp)} f_{G\mu\alpha}, \Pi_p) = \Tr(\chi_b^{G(\qp)} f_{G\mu\alpha}, \Pi_{0,p})$ for all contributing $\Pi$. We conclude that the trace in Equation~\eqref{eqn:Symp_expanded} is non-zero. This completes the proof.
\end{proof}

\begin{proof}[Proof of the existence of $\Pi_0$.] To find $\Pi_0$ we use the work of Shin~\cite{shinplancherel} (see in particular Thm~5.8 of [\textit{loc. cit.}] and its proof). However, we have to be a little bit careful, because Shin works under the condition that the center of $G$ is trivial. For most groups $G$ the existence of $\Pi_0$ follows from Shin's result: In case the center $Z$ of $G$ has $\uH^1(\Q, Z) = 1$ (and for all $v$, $\uH^1(\Q_v, Z) = 1$), to give an automorphic representation of $G/Z$ is to give an automorphic representation of $G$ whose central character is trivial. The center of the symplectic group of similitudes is $\Gm$, similarly the center of the odd spinor groups of similitudes is also a split torus. Thus for these groups the various cohomology groups $\uH^1$ vanish. However, for the unitary groups and the even spinor groups this is false for example. If $G$ has this property for its center, then we can apply Shin's result to $G/Z$ to find a suitable automorphic representation $\Pi_0'$ of $G/Z$. Lift it via $G(\A) \surjects G/Z(\A)$ to get the sought for automorphic representation $\Pi_0$ of $G(\A)$.

Now drop the assumptions on the $\uH^1$ of the center of $G$. We prove existence of $\Pi_0$. The argument and some of the text below is copied from Shin's article~\cite{shinplancherel}. We need to make some (very) minor adjustments to Shin's argument due to the presence of the center in our group $G$. We first simplify the notations somewhat. Write $\Sigma := S \cup \{p, p_1, p_2\}$. We write $\widehat U$ for the set of representations $\Pi$ of $G^*(\Q_\Sigma)$ such that
\begin{itemize}
\item $\Pi_p \in \iR(b)$;
\item $\Pi_{p_1}$ is an unramified twist of the Steinberg representation of $G(\Q_{p_1})$;
\item $\Pi_{p_2}$ lies in the inertial orbit of the representation $\Pi_{p_2}^\sharp$;
\item $\Pi_{S}$ is a twist of the Steinberg representation of $G^*(\Q_S)$.
\end{itemize}
The set $\widehat U$ has positive Plancherel density in the unitary dual $\widehat {G^*(\Q_\Sigma)}$. 

We claim (cf. Shin~\cite[thm 5.8]{shinplancherel}) that there exists infinitely many automorphic representations $\Pi$ of $G(\A)$ such that 
\begin{itemize}
\item $\Pi^{K^*} \neq 0$;
\item $\Pi_S \in \widehat U$;
\item $\Pi_\infty$ is $\xi$-cohomological for an irreducible representation $\xi$ of $G_\C$ which is sufficiently regular and has trivial central character. 
\end{itemize}
We prove the claim using Shin's notations/conventions (and we copy some of his text): Let $\xi_n$ be a sequence of irreducible representations of $G_\C$ with $\lim_{n \to \infty} \xi_n = \infty$, $\xi_n$ is sufficiently regular for all $n$, and the central character of $\xi_n$ is trivial for all $n$. Suppose the claim is false. Shrink $\widehat U$ so that it contains no automorphic representations but still has positive density. Then $m_{\textup{cusp}}(\Pi_S^0; \cchar_{K^{*S\infty}}, \xi_n) = 0$ for all $n \geq 1$. This implies that $\widehat \mu^{\textup{cusp}}_{\cchar_{K^{*S\infty}}, \xi_n}(\widehat U) = 0$. This contradicts the conclusion of Proposition~4.12 of [\textit{loc. cit}]. Unfortunately, Prop.~4.12 is proved under the assumption that the center of $G$ is trivial, so we need to show that a suitable conclusion is still satisfied in our setting. Prop.~4.12 refers to Theorem 4.11 of [\textit{loc. cit.}] and in the remark below the proof of Thm.~4.11, Shin explains that he needs the condition on the center to assure the vanishing of a certain limit $I_{2,n}$. This vanishing is in turn proved in Lemma 4.9 of [\textit{loc. cit.}]. Because we take $\xi_n$ to have trivial central character, it is clear that the limit in point (\textit{i}) of that lemma vanishes if $\gamma$ lies in the center (see also the remark 4.10 of [\textit{loc. cit.}]). Thus the conclusion of Prop.~4.12 holds also in our setting. The proof is complete.
\end{proof}

\section{Symplectic groups}
Let $G$ be a group of a Shimura datum satisfying the conditions from \S 1, and additionally the following: We have an exact sequence $[G_1 \injects G \surjects \Gm]$, where $G_1$ is an inner form of a symplectic group over a totally real extension $F^+$ of $\Q$. We call the morphism $G \surjects \Gm$ the factor of similitude. We show that $G$ satisfies Hypothesis~\ref{hypG} under these assumptions completing the proof of the nonemptiness for PEL varieties of type C.

Using results of Arthur~\cite{arthurbook} we show that Hypothesis~\ref{hypG} is true for $G$. Arthur recently proved his conjecture on the discrete spectrum for the group $G_1^*$. This is sufficient information for us. We write $G^*$ (resp. $G_1^*$) for the quasi-split inner form of $G$ (resp. $G_1$). The dual group $\widehat {G_1^*}$ is a product of $[F^+:\Q]$ copies of the group $\SO_{2g+1}(\C)$, with the absolute Galois group acting by permutation of the factors. We write ${}^LG_1^* := \SO_{2g+1}(\C)^{[F^+:\Q]} \rtimes \cG_\Q$. We recall briefly the result of Arthur for the group $G_1^*$. Let $\pi$ be a discrete automorphic representation of the group $G_1^*$. Arthur associates to $\pi$ a parameter $\psi \colon \cL_{\Q} \times \SU_2(\R) \to {}^{L} G_1^*$, where $\cL_\Q$ is the conjectural\footnote{We refer to \cite[\S 1.3]{arthurbook} for the explanation how these statement can be made independent of conjectures (at the cost of some technical burden). Roughly, one replaces the parameters $\psi|_{\cL_{\Q}}$ with essentially self-dual cuspidal automorphic representations of $\Gl_N$.} global Langlands group of $\Q$. These $L$-parameters $\psi$ are determined (up to conjugation) by the set of their local components $\psi_v \colon \cL_{\Q_v} \times \SU_2(\R) \to {}^{L} G_1^*$, where $v$ ranges over a (any) set of places \emph{with finite complement}, where $\cL_{\Q_v}$ is the local Langlands group. Using endoscopy Arthur parametrizes the elements of the packet $\Pi(\psi)$. Write, if $v$ is a $\Q$-place, $\cG_{\Q_v}$ for the absolute Galois group of $\Q_v$. For each place $v$ define $\widehat \cS_{\psi, v} := \Hom(\textup{Cent}_{\widehat G}(\psi) / \textup{Cent}_{\widehat G}(\psi)^\circ Z(\widehat G)^{\cG_{\Q_v}}, \C^\times)$. For each element of the packet $\pi_v \in \Pi(\psi_v)$, there is a character $\langle \cdot, \pi_v \rangle \in \widehat \cS_{\psi, v}$. In the tempered case, the characters in $\widehat \cS_{\psi}$ parametrize the elements of the local packets $\Pi(\psi_v)$. The global packet $\Pi(\psi)$ is the set of restricted tensor products $\pi = \bigotimes_v^\prime \pi_v$ such that
\begin{itemize}
\item $\pi_v \in \Pi(\psi_v)$;
\item The representation $\pi_v$ is unramified for nearly all $v$. 
\end{itemize}
Arthur defines the subset $\Pi(\psi)(\eps_{\psi}) \subset \Pi(\psi)$ consisting of those $\pi$ such that additionally:
\begin{itemize}
\item The product $\prod_v \langle \cdot, \pi_v \rangle$ is equal to a certain~\cite[p.~47]{arthurbook} explicit character $\eps_{\psi}$ depending only on the parameter $\psi$. 
\end{itemize}
Then he proves that $\Pi(\psi)(\eps_{\psi})$ is precisely the set of $\pi \in \Pi(\psi)$ occurring in the discrete spectrum of $G_1$. Arthur also computes the multiplicities (they turn out to be 1), but we will not need that fact.

\begin{theorem}[Conditional on the stabilization of
the trace formula]\label{hyptruegsp}
Hypothesis~\ref{hypG} is true for the group $G$.
\end{theorem}
\begin{remark}
The conditional character is due to the fact that Arthur's results are conditional.
\end{remark}
\begin{proof}[Proof of Theorem~\ref{hyptruegsp}.]
Let $(\Pi, \Pi_0)$ be a couple of cuspidal automorphic representations of the group $G^*$, satisfying the conditions of Hypothesis~\ref{hypG}. By assumption, the representation $\Pi_v$ is a twist of $\Pi_{0, v}$ at the places $v \in \Sigma := \{x_1, \ldots, x_r, p_1, p_2\}$. We must show that $\Pi_p$ is an unramified twist of $\Pi_{0, p}$. To simplify notations, we assume in the rest of the proof $\Pi_v \cong \Pi_{0, v}$. Furthermore, we only have to prove a property of automorphic representations of the quasi-split inner forms of $(G_1, G)$. Thus, to simplify notations, we can assume that $G_1$ and $G$ are already quasi-split themselves. 

By Proposition~\ref{st:extension} there exist cuspidal automorphic representations\footnote{Let us emphasize that the cuspidal automorphic $G_1(\A)$-subrepresentations $\pi$ of $\Pi$ are by no means unique or canonical.} $\pi, \pi_0$ of $G_1$ contained in $\Pi, \Pi_0$. Let $\psi$ and $\psi_0$ be the Arthur parameters for $\pi$ and $\pi_0$. The representations $\pi, \pi_0$ are automorphic, and thus unramified outside a finite set of places\footnote{Beware that the set of places where $\pi$ (resp. $\pi_0$) ramifies may very well be larger than the set of places where $\Pi$ (resp. $\Pi_0$) ramifies.}. 
We know that $\Pi$ and $\Pi_0$ are isomorphic outside a finite set, and furthermore, a local unramified representation of $G(\Q_v)$ (for some $v$) contains exactly one unramified representation of $G_1(\Q_v)$. Hence $\pi$ and $\pi_0$ are cuspidal representations of $G_1$ which are isomorphic outside a finite set of places. Consequently $\pi$ and $\pi_0$ have the same Arthur parameter: $\psi \colon \cL_\Q \times \SU_2(\R) \to {}^L G_1$. We have $\pi_p, \pi_{0, p} \in \Pi(\psi_p)$, and we want $\pi_p \cong \pi_{0, p}$. We look first at the infinite place. 
At this point we use the regularity condition on $\Pi_\infty$. This condition eliminates cohomology outside the middle dimension. In particular, our parameter $\psi$ is tempered and thus trivial on the $\SU_2(\R)$ factor. 
Locally at $p$ this means that the second $\SU_2(\R)$-factor of the local parameter 
$$
\psi_p \colon \underset {\cL_{\Q_p}}{\underbrace{\cW_{\qp} \times \SU_2(\R)}} \times \SU_2(\R) \to {}^L G_1,
$$
is trivial. Thus we have the more fine information that $\pi_p, \pi_{0, p}$ lie in the same packet defined by a tempered parameter. Now assume that the $p$-th local component $\Pi_{0, p}$ lies in the set $\iR(b)$. The representation $\pi_{0, p}$ is a subrepresentation of $\Pi_{0, p} = \Ind_{P(\qp)}^{G(\qp)}(\St_{M(\qp)}(\chi))$ for some standard parabolic subgroup $P$ of $G_{\qp}$.

We have the following lemma:

\begin{lemma}\label{surprise}
Let $M \subset G_{\qp}$ be a standard Levi-subgroup. Let $M_1 \subset M$ be the kernel of the morphism of similitudes. Let $\St_{M(\qp)}(\chi)$ be an unramified unitary twist of the Steinberg representation of $M(\qp)$.
\begin{itemize}
\item[(\textit{i})] 
The restriction of the representation $\St_{M(\qp)}$ to $M_1(\qp)$ is irreducible. 
\item[(\textit{ii})] 
Assume that the representation  $\Ind_{P(\qp)}^{G(\qp)}(\St_{M(\qp)}(\chi))$ is irreducible. The restriction of the induced representation $\Ind_{P(\qp)}^{G(\qp)}(\St_{M(\qp)}(\chi))$ to $G_1(\qp)$ is irreducible and isomorphic to $\Ind_{P_1(\qp)}^{G_1(\qp)}(\St_{M_1(\qp)}(\chi'))$, where $\chi'$ is the restriction of $\chi$ to $M_1(\qp)$.
\end{itemize}
\end{lemma}

(The lemma is proved at the end of this section.)
\smallskip

Consequently $\pi_{0, p}$ is an induction of the form $\Ind_{P_1(\qp)}^{G_1(\qp)}(\St_{M_1(\qp)}(\chi'))$. The Steinberg representation $\St_{M_1(\qp)}$ is the unique representation in its packet, and the parabolic induction of a tempered packet is again a packet. The Langlands packet of a parabolically induced representation $\Ind_{P}^G(\rho)$ consists precisely of all the irreducible representations that occur in parabolic inductions $\Ind_{P}^G(\rho')$, where $\rho'$ ranges over the representations lying in the same Langlands packet as $\rho$. By assumption on $\chi$, the induced representation $\pi_{0, p} = \Ind_{P_1(\qp)}^{G_1(\qp)}(\St_{M_1(\qp)}(\chi'))$ is irreducible. Thus the tempered packet in which $\pi_{0,p}$ occurs contains exactly one element. Since $\pi_p$ lies in the same packet as $\pi_{0, p}$, and this packet is a singleton, we get $\pi_p \cong \pi_{0, p}$. We have $\pi_p \subset \Pi_p$, and by applying Lemma~\ref{surprise} again,  $\Pi_p$ is isomorphic to $\Pi_{0, p}$ up to an unramified twist by some character $\chi$. Since the representation $\xi$ at infinity has trivial central character, this character $\chi$  must be of finite order. The proof is now complete. 
\end{proof}

\begin{proof}[Proof of Lemma~\ref{surprise} (i).]
In this proof we identify groups with their sets of $\qp$-rational points (we write $G, M, P, \ldots$, where we should have written $G(\qp), M(\qp), \ldots$). Let $B$ be the standard Borel subgroup of $M$, and write $B_1 := B \cap M_1$. We show that $M/B \isomto M_1/B_1$. We have the sequence $S_1 := [B_1 \injects B \surjects \qpm]$ mapping into the sequence $S_2 := [M_1 \injects M \rightarrow \qpm]$. Let us check that $B \rightarrow \qpm$ is surjective\footnote{The map $M \rightarrow \qpm$ is surjective as well, but we do not need that fact.}: We may decompose $B = TN$, $B_1 = T_1 N$, where $T_1$ the diagonal torus of $M_1$, $T$ the diagonal torus of $M$, and $N$ is the unipotent part (which is the same for $B$ and $B_1$). To check that $B \to \qpm$ is surjective, it suffices to check that $T \surjects \qpm$. To see that $T \surjects \qpm$, it suffices that $\uH^1(\qp, T_{1,\qp})$ vanishes. The torus $T_1$ is the restriction of scalars of a split torus, thus $\uH^1(\qp, T_{1,\qp})$ vanishes by Shapiro and Hilbert 90. Thus $B \rightarrow \qpm$ is surjective as claimed. By easy diagram chase we get from the morphism of sequences $S_1 \rightarrow S_2$ that $M/B = M_1/B_1$. Thus we have an $M_1$-equivariant isomorphism $C^\infty(M/B) \isomto C^\infty(M_1/B_1)$ between the locally constant complex valued functions on $M/B$ and $M_1/B_1$. 
Similarly $C^\infty(M/Q) \isomto C^\infty(M_1/Q_1)$ for any parabolic subgroup $Q \subset M$. Thus $C^\infty(M/B) \isomto C^\infty(M_1/B_1)$ induces a bijection from $\St_M$ to $\St_{M_1}$. Consequently, the restriction of $\St_M$ to $M_1$ equals $\St_{M_1}$.

Proof of (\textit{ii}). 
Consider the map $\Psi \colon \Ind_P^G(\St_M) \rightarrow \Ind_{P_1}^{G_1}(\St_M)$, defined by restricting a function $f \in\Ind_P^G(\St_M)$ to the subgroup $G_1 \subset G$. Clearly $\Psi$ is $G_1$-equivariant. We check by hand that $\Psi$ is bijective. First injective. Assume $f|_{G_1} = 0$. Let $g \in G$. Since $G/P = G_1/P_1$, we can write $g = p g_1$ with $p \in P_1$ and $g_1 \in G_1$. Then $f(g) = \sigma(p) f(g_1) = 0$. Therefore $f \equiv 0$, and $\Psi$ is injective. Proof that $\Psi$ is surjective: Let $f \in \Ind_{P_1}^{G_1}(\St_M)$, again because $G/P = G_1/P_1$, we may write for any $g \in G$, $g = pg_1$ with $p \in P_1$ and $g_1 \in G_1$. Define $\til f(g) := \sigma(p) f(g_1)$. 
Check that $\til f$ is welldefined: If $g = p'g_1' = p g_1$ are two decompositions, then $g_1' = (p')^{-1}p g_1$, and therefore $f(g_1') = \sigma((p')^{-1}p)f(g_1)$. Hence $\sigma(p')f(g_1') = \sigma(p) f(g_1)$, and $\til f$ is welldefined. Clearly $\Psi(\til f) = f$. Thus $\Psi$ is a bijection. By (\textit{i}) we have $\Ind_{P_1}^{G_1}(\St_{M}) \cong \Ind_{P_1}^{G_1}(\St_{M_1})$; compose with $\Psi$ to get $\Ind_P^G(\St_M) \cong \Ind_{P_1}^{G_1}(\St_{M_1})$. 
\end{proof}


\section{Odd Spin groups}

We prove non-emptiness of Newton strata for certain Shimura varieties associated to odd GSpin groups (so that the root system is of type B). 

Let us first introduce the Shimura varieties. Let $F^+$ be a totally real field, and let $(V, Q)$ be a quadratic space of $F^+$. We assume that $V$ is of odd dimension $n+2$, and that for every embedding $F^+ \hookrightarrow \R$ the associated real quadratic space $(V, Q)_\R$ is of signature $(n, 2)$. Let $C$ be the Clifford algebra over $F^+$ corresponding to $(V, Q)$, and write, for every $F^+$-algebra $R$, $C_R^+$ for the positive part of $C_R$. Let $G/\Q$ be the algebraic group such that for every $\Q$-algebra $R$, $\GSpin(V, Q)(R)$ is the set of $x \in (C_{R \otimes_\Q F^+}^+)^\times$ such that $x V_{R\otimes_\Q F^+} x^{-1} = V_{R \otimes_\Q F^+}$. Thus $G$ is the restriction of scalars to $\Q$ of a spinor group of similitudes over $F^+$. Let $X$ be the space of oriented negative definite $2$-planes in $V_{\R \otimes_\Q F^+}$. The couple $(G, X)$ is then a Shimura datum (see~\cite[proof of Lemma 3.2]{madapusipera}). The reflex field of $(G, X)$ is $\Q$. 

We reuse again the notations (and assumptions) that were introduced in \S 1 for general Shimura data.

\begin{theorem}[Conditional]\label{theorem_gspin}
Let $b \in B(G_{\qp}, \mu)$ be an admissible isocrystal. There exists a choice of Hecke operator $f^p$ and local system $\cL$ such that Formula~\eqref{equnow34} is non-zero.
\end{theorem} 
\begin{remark}
Conditional on the Rapoport-Langlands conjecture, the above theorem implies that the Newton strata are nonempty.
\end{remark}
\begin{remark}
Theorem~\ref{theorem_gspin} is conditional on the stabilization of the twisted trace formula (because we appeal to Arthur's description of the discrete spectrum).
\end{remark}
\begin{remark}
Due to problems with the center we show that the group $G$ satisfies only a slightly weaker version of Hypothesis~\ref{Ghyp}, as we only show the conclusion of Hyp.~\ref{Ghyp} holds for couples $(\Pi, \Pi_0)$ where the central character of $\Pi_0$ is trivial. We then show how the proof of Theorem~\ref{conditional_theorem} can be adjusted so that its conclusion is still valid under this weaker assumption on $G$. (See also the second paragraph of the proof below.) 
\end{remark}

\begin{proof}[Proof of Theorem~\ref{theorem_gspin}.]
Let $(\Pi, \Pi_0)$ be a couple of automorphic representations of $G$ as in~\ref{Ghyp}, and assume $\Pi_0 \in \iR(b)$. On top of the conditions in Hypothesis~\ref{Ghyp}, we assume also that the central character of $\Pi_0$ is trivial. We know that $\Pi^S \cong \Pi_0^S$. Thus the central character $\omega$ of $\Pi$ is trivial outside $S$. The infinite part of the central character is trivial as well (because $\xi$ has trivial central character). Thus $\omega$ factors over the group $Z(\R) Z(\Af^{S \cup \{p\}}) Z(\Q)$. By weak approximation, the group $Z(\Q)$ is dense in $Z(\A_{S \cup \{p\}})$, and thus $\omega$ is trivial. The center $Z$ of the group $G$ is equal to $\Res_{F^+/\Q} \Gm$ (see e.g.~\cite[Prop.~2.3]{MR2219256}; this uses that $n$ is odd). Therefore its first cohomology groups $\uH^1$ vanish over fields, and $Z(\A)G(\Q) \backslash G(\A) = O(\Q)\backslash O(\A)$, where $O := G/Z$. Thus $L^2(Z(\A)G(\Q) \backslash G(\A)) = L^2(O(\Q)\backslash O(\A))$. View $\Pi, \Pi_0$ as subspaces of the $L^2$-space of $O$. In odd dimension the special orthogonal group is adjoint. Therefore $O$ is the restriction of scalars from $F^+$ to $\Q$ of the quasi-split special orthogonal group $\SO_{n+2}$ over $F^+$. Thus Arthur's book~\cite{arthurbook} applies to $O$. By assumption we know that $\Pi_v$ and $\Pi_{0, v}$ coincide up to twist by certain characters $\chi_v$ if $v$ is one of the places $v \in \Sigma := \{x_1, \ldots, x_r, p_1, p_2, \infty\}$. For $v \notin \Sigma \cup \{p\}$ we have $\Pi_v \cong \Pi_{0, v}$ (also by assumption). Thus, the only place where we have a priori no information is the place $v = p$. Twist away the characters $\chi_v$ from the local representations $\Pi_v$, so that we have $\Pi^p \cong \Pi_0^p$. Let $\psi$ be the Arthur parameter for the representations $\Pi, \Pi_0$. We know that $\langle \cdot, \Pi_v \rangle = \langle \cdot, \Pi_{0, v} \rangle$ ($v \neq p$) and by Arthur, $\prod \langle \cdot, \Pi_v \rangle = \prod \langle \cdot, \Pi_{0, v} \rangle$ (products over \emph{all} $v$). Thus $\langle \cdot, \Pi_p \rangle = \langle \cdot, \Pi_{0, p} \rangle$. Since $\Pi_{0,\infty}$ and $\Pi_{\infty}$ are isomorphic, and both sufficiently regular by assumption, the parameter $\psi$ is tempered and trivial on the $\SU_2(\R)$-factor. Since $\Pi_0$ and $\Pi$ lie in the same tempered packet and $\langle \cdot, \Pi_p \rangle = \langle \cdot, \Pi_{0, p} \rangle$, we conclude (by Arthur) that $\Pi_p \cong \Pi_{0,p}$. This \emph{almost} proves Hypothesis~\ref{Ghyp} for the group $G$.

We said ``almost'' because we proved the hypothesis under an additional assumption on the central characters of $\Pi, \Pi_0$; we have to account for this assumption. Instead of trying to prove that Hyp.~\ref{Ghyp} is also true for representations without conditions on the central character, we alter the proof of Theorem~\ref{conditional_theorem} and show that we already have enough information to get the desired nonemptiness conclusion. The sequence $[\A^\times_{F^+} \injects G(\A) \surjects O(\A)]$ is exact. At the end of the proof of Theorem~\ref{conditional_theorem} where we prove the existence of $\Pi_0$, do not construct $\Pi_0$ as an automorphic representation of $G(\A)$, but instead as one of $O(\A)$. Thus, change the existence proof by replacing the role of $G(\A)$ by that of $O(\A)$. Then we have an automorphic representation $\Pi_0'$. Lift $\Pi'_0$ to an automorphic representation $\Pi_0$ of $G(\A)$ via the surjection $G(\A) \surjects O(\A)$. Then $\Pi_0$ has the conditions needed in the proof of Theorem~\ref{conditional_theorem} and $\Pi_0$ has furthermore trivial central character. Any representation $\Pi$ contributing to Equation~\eqref{eqn:Symp_expanded} has $\Pi^S \cong \Pi_0^S$ (where $S$ is some finite set of places), and therefore the central character of $\Pi$ is trivial outside $S$. Because the central character is also trivial at infinity, this implies that the central character of $\Pi$ itself is trivial. Thus, after these changements, we need in the proof of Theorem~\ref{conditional_theorem} the truth of Hypothesis~\ref{Ghyp} only for the couples $(\Pi, \Pi_0)$ where the central characters are trivial. The proof is complete. 
\end{proof}


\section{Appendix}
We prove two technical propositions, one concerning the restriction of automorphic representations to subgroups, and one proposition concerning certain Steinberg coefficients for groups with non-compact center.

\subsection{Restriction of automorphic representations}
Let $G$ be a connected reductive group over $\Q$ fitting in an exact sequence $[G_1 \injects G \surjects S]$, where $S$ is a torus of the form $\Res_{F/\Q} \Gm$, where $F/\Q$ is a finite \'etale $\Q$-algebra. We prove a technical proposition concerning the restriction of automorphic representations of $G$ to the subgroup $G_1 \subset G$. We write $Z$ (resp. $A$) for the center (resp. split center) of $G$. We write $Z_1, A_1$ for the corresponding objects of $G_1$. The mapping $G_1 \times Z \surjects G$ is surjective (over $\C$). Let $\Pi$ be a representation of $G(\A)$. View $\Pi$ by restriction via the mapping $G_1(\A) \times Z(\A) \to G(\A)$ as a representation $\Res(\Pi)$ of the group $G_1(\A) \times Z(\A)$.

\begin{proposition}\label{st:extension}
Let $\Pi$ be a cuspidal automorphic representation of $G(\A)$, then its restriction to the group $G_1(\A) \times Z(\A)$ contains a cuspidal automorphic representation of $G_1(\A) \times Z(\A)$. 
\end{proposition}
\begin{remark}
The proof we give here is copied from Clozel's article \cite[p.~137]{MR1114211}; cf. Labesse-Schwermer \cite[p.~391]{MR818358}. 
\end{remark}
\begin{proof}[Proof of Proposition~\ref{st:extension}] Define ${\mathbb G}_1$ to be the subset ${\mathbb G_1} := A(\A) G(\Q) G_1(\A) \subset G(\A)$. Then ${\mathbb G}_1$ is a subgroup because $G(\Q)$ normalizes $G_1(\A)$. Furthermore the subgroup ${\mathbb G}_1$ is closed in $G(\A)$, and we have $\left[ A(\A) G_1(\A) \right] \cap G(\Q) = A(\Q) G_1(\Q) \subset G(\A)$ (cf. Clozel [Lemme 5.8, \textit{loc. cit.}]). Let $\chi$ be the central character of $\Pi$; and let $\eps$ be the restriction of $\chi$ to $G_1(\A) \times A(\A)$. Let $\rho_0$ be the representation of $G_1(\A)$ on the space $L^2_0(G_1(\Q) \backslash G_1(\A), \eps)$ of cuspidal functions transforming under $G_1(\A)$ via $\eps$. We extend the representation $\rho_0$ to a representation of ${\mathbb G}_1$ by defining: $\rho_1(z \gamma x) f(y) = \chi(z) f(\gamma^{-1} y \gamma x)$, for $z \in A(\A), \gamma \in G(\Q), x \in G_1(\A), y \in G_1(\A)$. We do not copy the verification that this representation is well-defined [\textit{loc. cit}, 5.16]. Define the representation $\rho = \Ind_{ {\mathbb G}_1}^{G(\A)}(\rho_1)$ of $G(\A)$. A computation shows that $\rho$ is isomorphic to the representation of $G(\A)$ on the space $L^2_0(G(\Q) \backslash G(\A), \chi)$ of functions on $G(\Q)\backslash G(\A)$ transforming via $\chi$ under the action of $A(\A)$. Consequently, if $\Pi$ occurs in the representation $\Ind_{{\mathbb G}_1}^{G(\A)}(\rho_1)$, then its restriction to ${\mathbb G}_1$ will contain irreducible ${\mathbb G}_1$-subrepresentations of $\rho_1$. 
\end{proof}

\begin{remark}
The above proof has an interesting consequence. The argument shows that the $L_0^2$-space of $G_1$ carries a natural action of the group ${\mathbb G}_1$. Thus for any $g \in {\mathbb G}_1$ and any cuspidal automorphic representation $\pi$ of $G_1(\A)$, we can conjugate $\pi$ with $g$ to get a new cuspidal automorphic representation ${}^g \pi$ of $G_1(\A)$. (cf. \cite{MR818358}). 
\end{remark}

\subsection{Pseudocoefficients of the Steinberg representation}\label{appendix_section2}
We construct certain pseudo-coefficients $f_{x_i}$ of the Steinberg representation, in case the group has non-compact center. In the literature these coefficients are usually only constructed for groups under conditions on the center \cite[\S 3.4]{clozelchenevier}, such as the group be semi-simple, or with anisotropic center. We have neither of these conditions.

Let $F$ be a non-Archimedean local field, and let $G$ be the set of rational points of a connected reductive group over $F$. Let $H$ be the derived group of $G$. We write $G^*$ (resp. $H^*$) for the quasi-split inner form of $G$ (resp. $H$). We write $Z$ for the center of $G$. We write $H^*$ for the derived group of $G^*$ (then $H^*$ is the quasi-split inner form of $H$). The center $Z$ of $G$ is canonically isomorphic with the center of $G^*$ (and the same is true for the centers of $H$ and $H^*$). 

\begin{definition}\label{Average_over_Center}
Let $k$ be any smooth function on the group $H$. Let $O \subset Z$ be the maximal compact open subgroup of the center $Z$ of $G$. We now define a function $\widetilde k$ on the group $G$. Equip the finite group $H \cap Z$ with the counting measure $\dd t$. Consider the function $k^\star(g, z) := \int_{H\cap Z} (k \times \one_O)(gt, zt) \dd t$ on the group $H \times Z$. The function is $H\cap Z$-invariant, and thus defines a function on the subgroup $(H \times Z)/(H \cap Z) \subset G$. We extend this function by $0$ to obtain the function $\til k$ on the group $G$.
\end{definition}

The group $H^*$ is the derived group of $G^*$. By the fundamental lemma we may transfer smooth functions on the group $G$ to functions on the group $G^*$, and similarly functions from the group $H$ to functions on the group $H^*$. The formula in Definition~\ref{Average_over_Center} makes sense if we  replace $H$ by its quasi-split inner form; thus we also have a construction $k \mapsto \widetilde{k}$ for smooth functions on $H^*$. The construction in Definition~\ref{Average_over_Center} is compatible with transfer of functions, \ie the function $(\widetilde {k})^{G^*}$ on $G^*$ has the same stable orbital integrals as the function $\widetilde {\lhk k^{H^*} \rhk}$ for all $k \in \cH(H)$. 

We now take the function $k$ on $H$ to be a pseudocoefficient of the Steinberg representation, which exists because the center of $H$ is anisotropic. Define $f := \widetilde k$.

\begin{proposition}\label{labesseSteinberg}
Consider connected reductive group over a local field $G$ with quasi split inner form $G^*$. Assume that $G$ (and thus $G^*$) has anisotropic center. The transfer of a pseudocoefficient of the Steinberg representation is again a pseudocoefficient of the Steinberg representation. 
\end{proposition}
\begin{proof}
Labesse~\cite[Prop. 3.9.1, 3.9.2]{MR1695940}.
\end{proof}


\begin{definition}
Let $P_0$ be a minimal parabolic subgroup of the group $G$. Let $I$ be the space of locally constant functions on $G/P_0$. Then $I$ has an unique irreducible quotient $\St_G$, the \emph{Steinberg representation} of $G$.
\end{definition}

\begin{remark}
In case $G$ is anisotropic, the Steinberg representation $\St_G$ is the trivial representation of $G$. 
\end{remark}

We show that the function $f^{G^*}$ is (essentially) a pseudocoefficient of the Steinberg representation. 

\begin{definition}
Let $\chi$ be a character of the group $G$. The character $\chi$ induces a character $\li \chi$ of the cocenter $C$ of the group $G$. We call the character $\chi$ \emph{unramified} if $\li \chi$ is trivial on the maximal compact open subgroup $K_C$ of $C$.  
\end{definition}

The constructions in this section will be applied to a global setting, where the groups $G$ and $H$ arise as local components of global groups. Let us now assume additionally that this is the case: There exist groups $\underline G, \underline H$ defined over $\Q$, such that the fiber of $\underline G, \underline H$ at some $\Q$-place $x$ is equal to $G, H$. The quasi-split inner forms $G^*, H^*$ then arise as the $\Qx$-points of the quasi-split inner forms of $\underline G, \underline H$. 

\begin{proposition}\label{prop:stein_coefficients}
The function $(f)^{G^*}$ has the following two properties:
\begin{enumerate}
\item[(\textit{i})] For every unramified character $\chi$ of $G$:
$\Tr(f^{G^*}, \St_{G^*}(\chi))  \neq 0$. 
\item[(\textit{ii})] For every smooth irreducible representation $\Pi$ occurring as the $x$-component of a cuspidal automorphic representation the trace $\Tr(f^{G^*}, \Pi)$ is, up to the sign $\eps_{P_0}$, a non-negative real number.
\end{enumerate}
\end{proposition}
\begin{remark}
In our isolation argument we will correct the Steinberg coefficients by the sign $\eps_{P_0}$, so that the trace in (\textit{ii}) lies in $\R_{\geq 0}$. (One could also chose to not correct this sign; the importance in the global argument is that all the occurring traces have the \emph{same} sign, not so much that they are all non-negative).  
\end{remark}
\begin{proof}
We first verify (\textit{ii}). Let $\Pi$ be a smooth irreducible representation of the group $G^*$, let $\theta_{\Pi}$ be its character. We assume that $\Pi$ is the $x$-component of a cuspidal automorphic representation $\Pi$ of the group $G^*$. Let $\pi_1, \ldots, \pi_d$ be the irreducible $H^* Z$-subrepresentations of $\Pi$, and let $\theta_1, \ldots, \theta_d$ be their characters. We have $\theta_{\Pi} = \sum_{i=1}^d \theta_i$. Then (modulo a positive constant depending on $\dd t$):
\begin{equation}\label{eqn:Swap_Sum_Integral}
\Tr(f^{G^*}, \Pi) = \int_{H^* Z} f^{G^*}(g) \sum_{i=1}^d \theta_i(g) \dd g = \sum_{i=1}^d \int_{H^* Z} f^{G^*}(g) \theta_i(g) \dd g. 
\end{equation}
By the definition of the function $f^{G^*}$ the summand on the right hand side equals, up to some positive constant, the trace of the pseudocoefficient of the Steinberg representation on the group $H$ against $\pi_i$. Such a trace is non-zero only if $\pi_i$ is isomorphic to one of the representations $V_P$ defined by Borel-Wallach~\cite[6.2.14]{MR1721403}. We show that $\pi_i$ must be the Steinberg representation. The representation $\Pi$ occurs as the component at $x$ of a cuspidal automorphic representation. Therefore $\Pi$ is \emph{unitary}. Thus the representation $\pi_i$ is unitary as well. By [6.4, \textit{loc. cit.}] the only representations $V_P$ which are unitary, are the Steinberg representation and the trivial representation. Let us exclude the trivial representation. Because $G^*$ is quasi-split, the Steinberg representation is infinite dimensional. By Clifford theory, all the representations occurring in $\Pi$ are conjugate under elements of the group $G^*$. Consequently, if one of the occurring representations is finite dimensional, then they are all finite dimensional. This means that $\Pi$ is finite dimensional and thus the representation $\Pi$ is finite dimensional. Thus $\pi_i$ cannot be trivial. 

We now verify (\textit{i}). The function $f^{G^*}$ is supported on the inverse image of $K_C$ in $G$. Because $\chi$ is unramified it is constant on the support of $f^{G^*}$. Therefore we have $\Tr(f^{G^*}, \St_{G^*}(\chi)) = \Tr(f^{G^*}, \St_{G^*})$.  We verify that the trace  $\Tr(f^{G^*}, \St_{G^*})$ is non-zero. Let $P_{0, x}$ be a Borel subgroup of $G_{\Q}^*$ and let $P_{0, x}'$ be the pull back of $P_{0, x}$ to $H^*_{\Q}$. Let $I$ be the space of locally constant complex valued functions on $G^*/P_{0, x}$ and $I'$ be the same space, but then for the group $H^*$. We extend any function $h \in I'$ by $0$ and this gives us the composition of maps $I' \injects I \surjects \St_{G^*}$.  This composition is trivial on the subspaces $C^\infty(H^*/P) \subset I'$ for any proper parabolic subgroup $P$ of $H^*$ containing $P_{0, x}^*$. We obtain an $H^*$-injection $\St_{H^*} \injects \St_{G^*}$. It follows from Equation~\eqref{eqn:Swap_Sum_Integral} that $\Tr(f^{G^*}, \St_{G^*}) \neq 0$. 
\end{proof}

\subsection{A lemma on topological groups}\label{cliffordGR}
We prove some results regarding the restriction of representations. The statements can undoubtedly be found in the literature. We only prove the minimum of what we need (which is very little). Let $G$ be a topological group and $H \subset G$ an open normal subgroup with finite quotient. Let $\Pi$ be an irreducible representation of $G$ in a topological vector space. We prove\footnote{We borrowed this argument from the user ``vytas'' on mathoverflow.} some elementary results concerning the $H$-representations occurring in the restriction $\Pi|_H$ of $\Pi$ to $H$. Clearly $\Pi$ is a finitely generated $\C[G]$-module, and since the ring $\C[G]$ is finitely generated over $\C[H]$, $\Pi|_H$ is finitely generated as $\C[H]$-module. Then $\Pi$ has an irreducible quotient $\Pi|_H \surjects \pi$. Let $K$ be the kernel of $\Pi|_H \to \pi$. For every $g \in G$ the quotient $\Pi/gK$ is irreducible and isomorphic to $\pi^g$. The kernel of the natural map $\Pi \to \bigoplus_{g \in G/H} \Pi/gK$ is $G$-invariant, and hence $0$. Thus the $g(\pi)$ ($g \in G$) are the irreducible subrepresentations of $\Pi|_H$. Thus, for any two $\pi, \pi' \subset \Pi|_H$ there exists some $g \in G$ such that $\pi' \cong (\pi)^g$. If $g, g' \in G$ differ by an element of $H$, then $(\pi)^g \cong (\pi)^{g'}$ for all $\pi \subset \Pi|_H$. Since $H \subset G$ is open, this means that if $g$ and $g'$ are close enough, then automatically $(\pi)^g \cong (\pi)^{g'}$.

\bibliographystyle{plain}
\bibliography{grotebib}

\bigskip

\noindent Arno Kret, \\
\noindent Universit\'e Paris-Sud, UMR 8628, Math\'ematique, B\^atiment 425, \\
\noindent F-91405 Orsay Cedex, France

\end{document}